\documentclass[11pt,reqno]{article}
\usepackage[matrix]{xy}

\usepackage{amsmath, amssymb, amsthm}
\theoremstyle{definition}
\numberwithin{equation}{subsection}
\newtheorem*{Definition}{Definition}
\theoremstyle{plain}
\newtheorem*{definition}{Definition}
\newtheorem{theorem}{Theorem}

\newtheorem{prop}[theorem]{Proposition}
\newtheorem{lemma}[theorem]{Lemma}
\newtheorem{Lemma}[theorem]{Lemma}

\newtheorem{cor}[theorem]{Corollary}

\newcommand\Id{\text{Id}}

\DeclareMathOperator{\GL}{GL}
\DeclareMathOperator{\spn}{span}
\DeclareMathOperator{\inc}{inc}
\DeclareMathOperator{\proj}{proj}
\DeclareMathOperator{\Length}{Length}
 
\newcommand\barM{\overline M}

\newcommand\Haus{{\mathcal H}} 
\newcommand{\example}{\bigskip\noindent{\bf Example. }}
 
\newcommand{\Z}{{\mathbb Z}}

\newcommand{\R}{{\mathbb R}}
\newcommand{\N}{{\mathbb N}}
\newcommand{\C}{{\mathcal{C}}} 
\newcommand{\F}{{\mathcal{F}}} 
\newcommand{\G}{{\mathcal{G}}} 
\newcommand{\D}{{\mathcal{D}}} 
\newcommand{\A}{{\mathcal{A}}} 
\newcommand{\M}{{\mathcal{M}}}
\newcommand{\NBC}{{nilpotent-by-cyclic }}
\newcommand{\ABC}{{abelian-by-cyclic }} 
\newcommand\suchthat{\bigm|}
\newcommand{\n}{\mathfrak{n}}
\newcommand{\g}{\mathfrak{g}} 
\newcommand{\h}{\mathfrak{h}} 
\newcommand{\re}{\mathfrak{r}} 
\newcommand{\QIc}{{quasi-isometric }}
\newcommand{\qic}{{quasi-isometric}}
\newcommand{\QIy}{{quasi-isometrically }}
\newcommand{\qi}{{quasi-isometry}} 
\newcommand{\QI}{{quasi-isometry }}
 
\newcommand{\textmatrix}[4]{\bigl(
\begin{smallmatrix} #1 & #2 \\ #3 & #4 \end{smallmatrix} \bigr)}

\begin{document}
\title{The Large Scale Geometry of\\ Nilpotent-by-Cyclic Groups} 
\author{Ashley Reiter Ahlin\\2809B Hazelwood Drive, Nashville, TN 37212\\ ashleyahlin@yahoo.com }

\maketitle
\begin{abstract}

A nonpolycyclic \NBC group $\Gamma$ can be expressed as the HNN extension 
of a
finitely-generated nilpotent group $N$.  The first main result is that
\QIc \NBC groups are HNN extensions of \QIc nilpotent groups.  The
nonsurjective injection defining such an extension induces an injective
endomorphism $\phi$ of the Lie algebra $\g$ associated to the Lie group in
which $N$ is a lattice. A normal form for automorphisms of nilpotent Lie
algebras--\emph{permuted absolute Jordan form}-- is defined and 
conjectured to be a \QI invariant.  We show that
if $\phi$, $\theta$ are endomorphisms of lattices in a fixed Carnot
group $G$, and if the induced automorphisms of $\g$ have the same permuted
absolute Jordan form, then $\Gamma_\phi$, $\Gamma_\theta$ are \qic.  Two
\QI invariants are also found: \begin{itemize} \item The set of
``divergence rates'' of vertical flow lines, $\D_\phi$ \item The ``growth
spaces'' $\g_n \subset \g$ \end{itemize}

These do not establish that permuted absolute Jordan form is a \QI
invariant, although they are major steps toward that conjecture.

Furthermore, the \QIc rigidity of finitely-presented \NBC groups is
proven: any finitely-presented group \QIc to a nonpolycyclic \NBC group is
(virtually-nilpotent)-by-cyclic.

\end{abstract}

\tableofcontents
%\listoffigures
%\listoftables

\section{Introduction} 
\label{section:intro} 

The large-scale geometry of a group is captured by the notion of \qi.
A map $f$ between metric spaces $X$ and $Y$ is a {\em quasi-isometry} if
there are constants $K, C, C' \geq 0$ such that:
\begin{itemize}
\item for all $x, y \in X$, $$\frac{1}{K} d_X(x,y) - C \leq d_Y(f(x),f(y))
\leq K d_X(x,y) + C,$$ and \item the $C'$ neighborhood
of $f(X)$ is all of $Y$.  \end{itemize}

Given a finitely-generated group $G$, suppose that $S_1$ and $S_2$ are two
generating sets for $G$.  Let $d_i$ be the word metric on $G$ induced by
$S_i$. Then $(G, d_1)$ and $(G, d_2)$ are quasi-isometric.  Thus, the
quasi-isometry type of a group is independent of choice of generating set.
In 1980, Gromov established \cite{G} that the class of finitely-generated
nilpotent groups is \emph{\QIy rigid};  that is, a group is
quasi-isometric to a nilpotent group if and only if it is virtually
nilpotent.  One component of this proof was Wolf's result \cite{Wolf} that
nilpotent groups have polynomial growth, which is a \QI invariant and thus
a first step towards classifying nilpotent groups up to \qi.  
Furthermore, Bass \cite{Bass} found a formula for the precise degree of
polynomial growth of a nilpotent group, in terms of its lower central
series.  These results inaugurated the project of studying the large-scale
geometry of solvable groups via quasi-isometries.

One of the first major results in this area was due to Farb and Mosher,
who considered the solvable Baumslag-Solitar groups $BS(1,n)$ \cite{BS1}.  
They established that this class is \QIy rigid and gave a complete
classification up to \qi. These groups are the simplest examples of a rich
class of (typically) non-nilpotent solvable groups called \NBC groups. A
group $\Gamma$ is \emph{\NBC}if there is an exact sequence $$ 1
\rightarrow N\rightarrow \Gamma \rightarrow \Z \rightarrow 1 $$ where $N$
is a nilpotent group.

A finitely-generated \NBC
group $\Gamma$ can be thought of in a different way.  It is also an
ascending HNN extension of a nilpotent group $N$ defined by the 
endomorphism $\phi$.  If $N$
has presentation $N = \langle R \suchthat S \rangle$, then $\Gamma$ is
given by the presentation: $$ \Gamma = \langle R, t \suchthat S,
tnt^{-1}=\phi(n), \mbox{ for }n \in N\rangle.$$ By a result of Bieri and Strebel \cite{BS}, if $\Gamma$ is
finitely presented, then it can be represented as the ascending HNN
extension of a (possibly different) nilpotent group which is finitely
generated. The group $\Gamma$ is nonpolycyclic if $\phi$ is not
surjective.

The case when $N$ is abelian, i.e., the \emph{abelian-by-cyclic} groups, 
was considered by Farb and Mosher \cite{ABC} in 1999.  Defined
like \NBC groups, finitely-generated \ABC groups are the HNN extensions of
finitely-generated abelian groups.  The action of the nonsurjective 
injection can be
specified by an $n \times n$ integer matrix $M$.  The \emph{absolute
Jordan form} of a matrix is obtained from its Jordan form by replacing
each diagonal entry with its absolute value. They showed that two
finitely-generated, nonpolycyclic \ABC groups are \QIc if and only if the
matrices defining them have integral powers with the same absolute Jordan
form.  Furthermore, they showed that this class of groups is \QIy rigid.  
That is, any group which is \QIc to a group in this class has a quotient
by a finite subgroup which is (virtually) one itself.  The questions of
classification and rigidity of \NBC groups were later posed by Farb and
Mosher (\cite{FMS}, Problem 3 and Question 2, respectively).

\subsection{Statement of Results}\label{subsection:results} 

The first main result of this paper is:

\begin{theorem}[Rigidity] \label{theorem:rigidity}
Let $\Gamma = \Gamma_{N,\phi}$ be a finitely-presented nonpolycyclic \NBC
group.  If $G$ is a finitely-generated group quasi-isometric to $\Gamma$,
then $G$ is the ascending HNN extension of a virtually nilpotent group.  
\end{theorem}

To begin classifying \NBC groups up to \qi, we first show that two \QIc
\NBC groups must be based on nilpotent groups which are themselves
quasi-isometric.  The classification of nilpotent groups up to \QI is
still a major open question.  For the purposes of classifying groups up to
\qi, we will restrict our attention to lattices in \emph{Carnot groups},
which are nilpotent Lie groups with a particularly nice nilpotent grading.  
Lattices in Carnot groups are \QIc if and only if they are lattices in the
same Carnot group.  Thus, for the remainder of the work on classification,
we make the slightly stronger assumption that the base nilpotent groups
are lattices in the same Carnot group.  We also restrict our attention to
\NBC groups defined by endomorphisms which act without unipotent part.
(See Subsection \ref{subsection:endo} for a precise definition.)

Every finitely-generated nilpotent group has a finite-index torsion-free
subgroup (\cite{Baum}, Theorem 2.1, citing Hirsch), to which it must
therefore be quasi-isometric.  Also, by a result of Malcev (\cite{Mal}; 
see also \cite{R}, Theorem 2.18)
every finitely-generated torsion-free nilpotent group $N$ is a lattice in
a connected, simply-connected nilpotent Lie group $G$.  Malcev also shows
that any (injective) endomorphism of the discrete nilpotent group $N$
extends to an (injective) endomorphism of the Lie group $G$.  Given such
an endomorphism $\phi$, consider the induced linear map $\phi^*$ which
acts on the Lie algebra $\g$. We have restricted our attention to \NBC
groups based on nilpotent groups which are lattices in the same Carnot
group $G$.  Thus, we are comparing maps of a fixed Lie algebra $\g$.

Such a map can be represented by some matrix $M \in \GL_n({\mathbb C})$
which is in Jordan form. If all the eigenvalues of $M$ are real, then this
matrix represents the map with respect to some basis $\{e_i\}$.  (If $M$
has complex eigenvectors, we carry out a similar procedure to what
follows.) Associated to the basis $\{e_i\}$ is a weight vector $w \in
\N^n$ which specifies the component of the nilpotent grading in which each
vector appears.  Given this weight data, we permute the basis vectors to
obtain a matrix in \emph{permuted Jordan form}.  From this form, we obtain
the \emph{permuted absolute Jordan form} of $M$ by replacing each diagonal
entry with its absolute value.  (See Subsection \ref{subsection:pajf} for
precise definitions.)

This form is a classifier of \NBC groups; that is, groups which have the 
same permuted absolute Jordan form are quasi-isometric.

\begin{theorem}[Permuted Absolute Jordan Form] \label{theorem:pajf} Let
$N_1$ and $N_2$ be two lattices in the same Carnot group $G$.  Let
$\phi_1$ and $\phi_2$ be injective, nonsurjective endomorphisms of $N_1$
and $N_2$ respectively, each acting without unipotent part.  Suppose that
there are integers $r_1, r_2$ such that $M_1^{r_1}$ and $M_2^{r_2}$ have
the same permuted absolute Jordan form.  Then $\Gamma_{N_1,\phi_1}$ and
$\Gamma_{N_2,\phi_2}$ are \qic.  \end{theorem}

In the case of \ABC groups, absolute Jordan form is a \QI invariant.  
The absolute Jordan form captures very specific information about a
matrix.  It is determined by the absolute values of the eigenvalues and
the dimensions of the corresponding root spaces.  This information
determines \QI type by identifying the rates at which vectors grow when
repeatedly multiplied by the matrix $M$.  A pure eigenvector grows as
$\lambda^t$, where $\lambda$ is the corresponding eigenvalue.  A vector in
the $\lambda$-root space grows as $t^n \cdot \lambda^t$ for some $n \in
\N$.  Farb and Mosher proved that these growth rates are \QI invariants of the \ABC
group and that these rates uniquely determine the absolute Jordan form of
$M$.

In the case of \NBC groups, similarly defined divergence rates are still
\QI invariants.  For each $x \in N$, consider the function $f_x(t) = 
d(0, \phi^t(x))$.  After defining a suitable equivalence relation for 
functions from $\R$ to $\R$, we consider the set of all divergence rates, 
up to this equivalence relation:
$$\D_\phi = \{[f_x(t)] \suchthat x \in N\}.$$
See Subsection \ref{subsection:divrates} for precise definitions.

\begin{theorem}[Divergence Rates are a Quasi-isometry Invariant]
\label{theorem:invariant1} Let $N_1$ and $N_2$ be lattices in the same
Carnot group $G$.  Let $\phi_1$ and $\phi_2$ be injective, nonsurjective
endomorphisms of $N_1$ and $N_2$ respectively, each acting without
unipotent part.  If $\Gamma_{N_1,\phi_1}$ and $\Gamma_{N_2,\phi_2}$ are
\QIc then the multisets of divergence rates $\D_{\phi_1}$ and
$\D_{\phi_2}$ are equal. \end{theorem}

Because the nilpotent group has interesting geometry of its own, the set
of divergence rates is not sufficient to specify the absolute Jordan
form.  We use the permuted absolute Jordan form in order to keep some
information about the geometry of the nilpotent group, and this further
information ensures that the divergence rates can be calcuated from the
permuted absolute Jordan form.  Nonetheless, the permuted absolute Jordan
form is not uniquely determined by the set of divergence rates, so, unlike
the \ABC case, this does not establish permuted absolute Jordan form as a
\QI invariant.

The set of divergence rates fails to determine the permuted absolute 
Jordan form in part because it fails to capture which divergence rates 
arise from points in the various levels of the nilpotent group.  Part of 
this data is found by considering the Lie subalgebras associated to 
various growth rates.  We consider the \emph{growth spaces}:
$$\g_\lambda = \{v \in \g \suchthat ||M^tv|| \preceq \lambda^t t^k\mbox{ 
for some } k \in \N \}.$$
(See Subsection \ref{subsection:growth} for precise definitions.)
Pansu has shown \cite{P} that quasi-isometric Carnot groups have isomorphic Lie 
algebras.
This implies that the isomorphism type of the growth spaces is a \QI 
invariant of the \NBC groups. 

\begin{theorem}[Growth Spaces are a Quasi-isometry Invariant]
\label{theorem:invariant2} Let $N_1$ and $N_2$ be lattices in the same
Carnot group $G$.  Let $\phi_1$ and $\phi_2$ be injective, nonsurjective
endomorphisms of $N_1$ and $N_2$ respectively, each acting without
unipotent part.  If $\Gamma_{N_1,\phi_1}$ and $\Gamma_{N_2,\phi_2}$ are
\QIc then each growth space $\g_{\lambda_1}$ of $\Gamma_{N_1,\phi_1}$
is isomorphic (as a Lie algebra) to some growth space $\g_{\lambda_2}$ of
$\Gamma_{N_2,\phi_2}$. \end{theorem}

The remainder of this section will contain an outline of the
classification results. The second section focuses on a single \NBC group,
first describing precisely how the permuted absolute Jordan form is
computed and then showing how that data determines the geometry of a
particular associated Lie group. In Section \ref{section:coarsetop}, a
geometric model space for the \NBC groups is described, followed by proofs
of Theorems \ref{theorem:pajf} - \ref{theorem:invariant2}.  Some
low-dimensional examples of \NBC groups are described in Section
\ref{section:examples}. A proof of the rigidity result follows in the last
section.

\subsection{Outline of the Classification} 

{\bf Step 1:} (Subsection \ref{subsection:model}) Given a \NBC group $\Gamma =
N_\phi$, where $N$ is a finitely-generated nilpotent group and $\phi$ is 
an injective
endomorphism of $N$, we construct a geometric model space $X$ which is
\QIc to $\Gamma$. Topologically $X = \R^n \times T$, where $T$ is the
Bass-Serre tree for $\Gamma$.  For each $x \in T$, the
horizontal slice $\R^n \times \{x\}$ has the metric given by the pullback
of the metric on $\R^n \times \{0\}$ via $(\phi^t)^*$, where $t =\mbox{
height of }x$.  This construction follows the same process as that used in
\cite{ABC}, but the metric on each slice is the non-isotropic geometry of
nilpotent groups which is described in Subsection \ref{subsection:kar}.

\bigskip
\noindent
{\bf Step 2:} (Subsection \ref{subsection:reducetoajf}) 
We extend the results used by Farb and Mosher in \cite{ABC} to show that 
the group defined by a given matrix is \QIc to the group defined by its 
absolute Jordan form.  This requires showing that the conjugating matrices 
preserve the nilpotent structure of the base group.

\bigskip
\noindent
{\bf Step 3:} (Subsection \ref{subsection:coarsetop})
Suppose $f$ is a quasi-isometry between two \NBC groups $\Gamma =
\Gamma_{N,\phi}$ and $\Gamma' =\Gamma_{N',\phi'}$.  With the additional 
condition that $\phi$ and $\phi'$ are not surjective, we apply the same
coarse topology as in \cite{ABC} to get a \QI between the Lie groups which
are hyperplanes in the model spaces: $ G = N \rtimes_\phi \R$ and $ G' =
N' \rtimes_{\phi'} \R$. This result depends upon the groups $\Gamma$ and
$\Gamma'$ being nonpolycyclic, which holds when the Bass-Serre trees
$T_\phi$, $T_{\phi'}$ have valence $v \geq 3$.  These hyperplanes
correspond to $\R^n \times l$ for some directed line $l \subset T$.  As a
result, we show that if two \NBC groups are quasi-isometric then their
base nilpotent groups must also be quasi-isometric.  Having reduced this
part of the problem to the (unsolved) quasi-isometric classification of
nilpotent groups, we thenceforth restrict our attention to \NBC groups
which are HNN extensions of two nilpotent groups which are lattices in the
same nilpotent Carnot group.

\bigskip
\noindent
{\bf Step 4:} (Subsections \ref{subsection:flow} - \ref{subsection:time2})
The \QI between the hyperplanes has even more structure. Under the
additional condition that $\phi$ and $\phi'$ act without unipotent part,
we show that vertical flow lines are coarsely preserved. For any \NBC
group, we calculate the divergence rates for pairs of vertical flow lines
and show that the set of such rates is a finite set which is determined by the
permuted absolute Jordan form (Subsections
\ref{subsection:order} - \ref{subsection:divrates}).
By considering the form of such divergence rates, we show that the time
change function induced by the \QI is linear, so the set of rates at
which vertical lines diverge from the flow line at the origin is a \QI
invariant, up to rescaling all the rates by a single power.

\bigskip
\noindent
{\bf Step 5:} (Subsection \ref{subsection:growth}) 
We use the results of Section \ref{section:geom} on divergence rates to 
show that the growth spaces are Lie subalgebras and then apply Pansu's 
theorem to establish Theorem \ref{theorem:invariant2}.

\section{Calculating the Quasi-isometry Invariants}
\label{section:ajf}\label{section:geom}
%Section 2

The main result of this section is contained in Theorem
\ref{theorem:divrates}, which establishes the set of ``divergence rates''
in Lie groups associated to \NBC groups.  These rates are one of the two
\QI invariants established in Section \ref{section:coarsetop}.  We also
establish, in Subsection \ref{subsection:pajf}, that two groups with the same
permuted absolute Jordan form are \qic.

In the early subsections of this section, we present some tools which will be
needed for the proof in Subsection \ref{subsection:divrates}.  In Subsection
\ref{subsection:kar}, we describe the geometry of left-invariant metrics on
nilpotent groups.  We establish some characteristics of endomorphisms of
torsion-free nilpotent groups in Subsection \ref{subsection:endo}, and describe
a needed assumption on the endomorphism in Subsection \ref{subsection:unip}.  
The definition of permuted absolute
Jordan form is found in Subsection \ref{subsection:pajf}. Subsection
\ref{subsection:linalg} contains a key result from linear algebra which
relates the structure of a nilpotent Lie group to the structure of any
endomorphism of its Lie algebra.  Subsection \ref{subsection:opjs} describes and 
strengthens a result on the relationship between a matrix and its absolute 
Jordan form, in preparation for the proof of Theorem \ref{theorem:pajf} in 
the following subsection. In Subsection \ref{subsection:order}, we
introduce an equivalence relation which allows us to distinguish between
different divergence rates.  Subsection \ref{subsection:vectorgrowth} expands
upon results of \cite{ABC} to describe the growth of a vector under
repeated application of a linear map.  Finally, all these pieces are used
in the final subsection to establish the set of divergence rates for a \NBC
group.

\subsection{Geometry of Nilpotent Groups} \label{subsection:kar}

We consider the class of left-invariant Riemannian metrics on
connected, simply-connected nilpotent Lie groups.  As for discrete 
groups, any two left-invariant 
metrics on the same Lie group are quasi-isometric, so the quasi-isometric 
classification of groups, independent of the particular metric, is a well-defined problem. 

Such groups can be can be globally coordinatized by $\R^n$ such
that balls centered at the origin are comparable to ellipsoids with axes 
of length which is polynomial in the radius,
relative to the coordinates of $\R^n$.

Nilpotent groups admit both Riemannian metrics and so-called
Carnot-Caratheodory (CC) metrics.  Typical CC metrics are non-Riemannian
on the infinitesimal scale and are non-isotropic on the large scale. The
Riemannian metrics of interest to us are trivially CC metrics.  They are
infinitesimally Euclidean, unlike typical CC metrics.  However, on the
large-scale, Riemannian metrics and CC metrics are similarly
non-isotropic.  This non-isotropic nature is revealed by the description
of balls in this metric, which is found below in Theorem 
\ref{theorem:kar}.

Given a finitely-generated nilpotent group $N$, we define subgroups $\gamma_1(N)=N$,
$\gamma_2(N)=[\gamma_1(N), \gamma_1(N)]$, and inductively,
$\gamma_{i+1}(N)=[\gamma_1(N), \gamma_i(N)]$.  Then $$N = \gamma_1(N)
\supset \cdots \supset \gamma_{c+1}(N) = 1$$ is the lower central series
of $N$. We will refer to this filtration on $N$ as the \emph{nilpotent
grading}. Define: $$d_i = \dim (\gamma_i/\gamma_{i+1}).$$

We choose a basis for the Lie algebra which respects this grading:

\begin{definition}[Triangular basis]  
Suppose $\{e_1, \ldots e_n\}$ is a basis for the nilpotent Lie algebra
$\n$ such that $[e_i, e_j] = \sum_k\alpha_{ijk}e_k$.  The basis is 
\emph{triangular} if $\alpha_{ijk} = 0 $ when $k \leq \max(i,j)$.  The 
constants  $\alpha_{ijk}$ are called the \emph{structure constants} for 
the group.
\end{definition}  

\example The vectors $\{X, Y, Z\}$ form a triangular basis for the Lie
algebra of the Heisenberg group, because $[X, Y] = Z$, $[X, Z] = 0$, and
$[Y, Z] = 0$.  The basis $\{X+Z, Y, Z\}$ is also triangular, but neither
$\{X, Y, X + Z \}$, nor the reordered basis $\{Z, Y, X\}$ is.
\bigskip

Each vector $v \in \n$ is assigned a weight $w(v)$ which specifies the last component of the nilpotent grading which contains $v$.  For basis vectors $e_k$:
$$w_k = w(e_k)) = \max \{i \suchthat e_k \in \gamma_i(\g) \}.$$
For a given choice of (ordered) basis, we will call the associated $n$-tuple $(w_1, \ldots, w_n)$ the \emph{weight vector} associated to the basis.

\example
In the Heisenberg group, given in $\{X, Y, Z\}$ coordinates, 
the weights are: $w(X) = w(Y) = 1; w(Z) = 2$.  That is, the weight vector 
is $(1, 1, 2)$. 
\bigskip

In this notation, both Gromov (\cite{Gballbox}, for Carnot-Caratheodory 
spaces)  
and Karidi have shown (\cite{K}, Theorem 4.2):

\begin{theorem}[Ball-Box Comparison Theorem]\label{theorem:kar} Let $N$ be 
a connected, simply-connected, real, nilpotent
Lie group of dimension $n$ with Lie algebra $\n$, and let $\{e_1, \ldots,
e_n\}$ be a triangular basis of $\n$ with associated weight vector
$w = (w_1, \ldots, w_n)$.  Specifying $\{e_1, \ldots, e_n\}$ as an orthonormal basis defines an inner product at the origin.  Applying left-invariance, this defines a metric on the group $N$.  Let $B(r)$ be the ball centered at the 
origin in $N$ with radius $r > 1$ in this metric.

Then there exists a constant $a >1$ (which depends on the group $N$,
but not on $r$) such that 
$$\{|x_i|\leq (r/a)^{w_i} \suchthat i = 1, \ldots, n\}
\subset B(r) \subset \{|x_i|\leq (ar)^{w_i}\suchthat i = 1, \ldots, n\}.$$
\end{theorem}

\begin{definition}
Two functions $a, b\colon X \rightarrow \R$ are \emph{comparable}, denoted
$a(x) \sim b(x)$ if there exists $K > 0$ such that for all $x \in X$
$$\frac{1}{K}b(x)  < a(x) < K b(x).$$
\end{definition}

\example
In the Heisenberg group, given in $\{X, Y, Z\}$ coordinates, 
the ball of radius $r$ is comparable to the box of the form $[-r, r] \times 
[-r, r] \times [-r^2, r^2]$.
\bigskip

Karidi's result on balls can be rephrased as follows to describe the
distance between points.  

\begin{cor}[Distances in Nilpotent Groups]\label{cor:dist1}
For a nilpotent Lie group $N$, with weights $w_i$ as in Theorem 
\ref{theorem:kar}, the left-invariant metric:
$$||(x_1, \ldots, x_n)|| \sim \max_i\{|x_i|^{\tfrac{1}{w_i}}\}.$$
\end{cor}

\begin{proof}
For $ x = (x_1, \ldots, x_n)$, let $d= \max_i\{|x_i|^{\tfrac{1}{w_i}}\}$.  
Then each $|x_i|^{\tfrac{1}{w_i}} \leq d$, so $|x_i| \leq d^{w_i} $.  By
Theorem \ref{theorem:kar}, the box $\{|x_i| \leq d^{w_i}\} \subset 
B(ad)$, so $x \in B(ad)$.  Conversely, there is some $i$ such that 
$|x_i|^{\tfrac{1}{w_i}} = d$, so  $|x_i| = d^{w_i} $.  Again, 
Theorem \ref{theorem:kar} implies that for all $\epsilon > 0$, $x \not\in
\{|x_i| \leq (d- \epsilon)^{w_i}\} $, which contains 
$B(\frac{d-\epsilon}{a})$.  Thus, $\frac{d-\epsilon}{a} < ||x|| < ad$ for 
all $\epsilon > 0$, so $\frac{d}{a} \leq ||x|| < ad$.
$$||(x_1, \ldots, x_n)|| \sim \max_i\{(x_i)^{\tfrac{1}{w_i}}\}.$$
\end{proof}

At times, it will be convenient to use the following characterization of 
this distance metric:

\begin{cor}[Distances in Nilpotent Groups, Part 2]\label{cor:dist2}
Given a triangular basis $\{e_i\}$ with weights $\{w_i\}$, consider 
$$V_j = \spn\{e_i \suchthat w_i = j\}.$$
Then for some $k$, $G \cong R^n = V_1 \oplus \cdots \oplus V_k$.
Given $x \in G$, we express $x$ uniquely as $x = (x_1, \ldots, x_k)$, 
with each $x_i \in V_i$.  Then,
$$||x|| \sim \max_i\{\sqrt[i]{|x_i|}\}.$$
\end{cor}

\begin{proof}
This follows immediately from the definitions and Corollary 
\ref{cor:dist1}.
\end{proof}

\example In the Heisenberg group, the distance from the origin to the
point $(a, b, c)$ (which represents the group element $x^a y^b z^c$) is
comparable to the function $\max\{a, b, \sqrt{c}\}$. 
\bigskip

The following four-step nilpotent group will be used as the basis 
for examples of \NBC groups in Subsection \ref{subsection:example2}.

\example
Define
\begin{equation}\label{eqn:4step}
\begin{split}
G = \langle x,y,z,a,b,c,p, q,r,s,t\suchthat& [x,y]=z, [a,b]=c, [z,c]=t,\\
&[p,q]=r, [p,r]=s,[q,r]=s,\\
&[r,r]=t, [p,s]=t, [q,s]=t \rangle.
\end{split}
\end{equation}

Notice that $(x,y,z)$, $(a,b,c)$, $(p,q,r)$ and $(z,c,t)$ are each
isomorphic to the Heisenberg group $H$.  This nilpotent group has:
\begin{equation*}
\begin{split}
V_1 &= \langle x,y, a,b, p, q\rangle\\
V_2 &= \langle z,c, r\rangle\\
V_3 &= \langle s\rangle\\
V_4 &= \langle t\rangle\\
\end{split}
\end{equation*}
Thus,
$$d(0, (x,y,z,a,b,c,t)) \simeq \max \{x,y,a,b,p,q, \sqrt{z}, \sqrt{c},
\sqrt{r}, \sqrt[3]{s}, \sqrt[4]{t}\}. $$
\bigskip

\subsection{Some Characteristics of Endomorphisms of Nilpotent Groups}
\label{subsection:endo} 
Recall that a \NBC group can be expressed as the HNN
extension of a finitely-generated nilpotent group $N$ by an injective 
endomorphism $\phi$.  
We will now consider the endomorphisms of nilpotent groups which can
define such an HNN extension.  Their classification is an open question
and is not considered here.  However, we will establish a few important
facts about such endomorphisms.

Malcev has shown (\cite{Mal}, Theorem 5) that an automorphism of a
discrete nilpotent group $\Gamma$ can be extended to an automorphism of
any nilpotent Lie group in which $\Gamma$ is a lattice.  This result has
been broadly generalized, for example by Raghunathan (\cite{R}, Theorem
2.11, p. 33).  We will use the following statement:

\begin{theorem}\label{theorem:malcev} Given an injective endomorphism 
$\phi$ of the discrete nilpotent group $\Gamma$, we can extend to an 
injective endomorphism $\hat\phi$
of the nilpotent Lie group $G$ in which $\Gamma$ is a lattice. As always,
this extends to an automorphism of the Lie algebra $\g$. \end{theorem}

The algebraic structure of a nilpotent group places significant
restrictions on the structure of its endomophisms.  The restrictions are
most easily described for a special class of nilpotent group known as
\emph{Carnot groups}, defined below.  As it turns out, these groups
are the ones which admit the largest classes of endomporhisms.

\begin{definition} A connected, simply-connected nilpotent Lie group $G$ is 
a \emph{Carnot group} if its Lie algebra $\g$ admits a grading 
$\g =\bigoplus_{j=1}^r V_j$ such that:
\begin{itemize}
\item $[V_1, V_j]=V_{j+1}$ for $1 \leq j < r$, and $[V_1, V_r]=\{0\}$, and
\item $V_1$ generates all of $\g$ via $[\cdot, \cdot]$.
\end{itemize}
\end{definition}

Any nilpotent Lie group has a naturally associated Carnot group, which is 
obtained by redefining each commutator, omitting any term which is in a 
higher grade than permitted by the Carnot definition.

\begin{definition}
Given a linear map $\phi$ of a vector space $V$, and a grading $V = 
\bigoplus_{j=i}^r V_j$, we say that $\phi$ \emph{weakly preserves} the 
grading of $V$ if $\phi(\bigoplus_{j=i}^r V_j) \subseteq (
\bigoplus_{j=i}^r V_j)$ for each $i=1, \ldots, r$.
\end{definition}

\begin{lemma} \label{theorem:endo} Let $G$ be a Carnot group, and let
$\phi$ be an injective endomorphism of the Lie algebra $\g$.  Then, $\phi$
weakly preserves the Carnot grading $\g = \bigoplus_{j=1}^r V_j$.  
\end{lemma}

\begin{proof} 
Proof is by induction.  The base case, that $\phi (\bigoplus_{j=1}^r V_j)
\subseteq \bigoplus_{j=1}^r V_j$ is trivially satisfied.

Now, assume that $\phi(\bigoplus_{j=i}^r V_j) \subseteq  (
\bigoplus_{j=i}^r V_j)$ for each $i=1, \ldots, k$.  
Each $v \in \bigoplus_{j=k+1}^r V_j$ can be expressed as a 
commutator $v = [x,y]$ for some $x \in \R^n, y \in \bigoplus_{j=k}^r V_j$.
Then, $\phi(v) = [\phi(x), \phi(y)]$, with
$\phi(x) \in \g$ and $\phi(y) \in \bigoplus_{j=k}^r V_j$, so $\phi(v) \in 
\bigoplus_{j=k+1}^r V_j$.
\end{proof}

The structure of the nilpotent group does give even more restrictions on
the endomorphism.  In particular, for a Carnot group, the action of $\phi$
on the base level, $V_1$, determines completely the action on the rest of
the group.

\example
An endomorphism of the Heisenberg group can be described 
via the $3 \times 3$ matrix which represents the endomorphism for the basis $\{x, y, z\}$.
Given the first two columns (representing the effect of $\phi$ on $x$ and 
$y$) the third column (representing the effect of $\phi$ on $z$) can be 
calculated explicitly, since $\phi(z) = [\phi(x), \phi(y)]$.  If $\phi(x) = x^ay^cz^e$ and 
$\phi(y) = x^by^dz^f$, then $\phi(z) = z^{det}$, where det $=\left|
\begin{array}{cc}
a&b\\
c&d\\
\end{array}
\right|$, the determinant of the matrix representing $\phi|_{<x,y>}$.  So, a generic matrix representing an endomorphism of the Heisenberg group has the form:
$\left[
\begin{array}{ccc}
a&b&0\\
c&d&0\\
e&f&\mbox{det}\\
\end{array}
\right].$
\bigskip

\subsection{Unipotent-Free Endomorphisms}\label{subsection:unip}

\begin{definition} An endomorphism $\phi$ of a nilpotent group $G$ is
called \emph{unipotent-free} if the matrix representing $\phi$ has no
unipotent part, that is, no Jordan block with eigenvalue on the unit
circle. \end{definition}

We now make a:\\
{\bf Standing Assumption:} The endomorphism $\phi$ of $G$ is 
unipotent-free.

This assumption was not necessary in the \ABC case.  The extension of an
abelian group by unipotent matrix $U$ is a nilpotent group of the form
$\Z^k \rtimes_U \Z$.  As shown by Bridson and Gersten \cite{BG}, Pansu's
invariant \cite{P} shows that the sizes of the unipotent blocks in $U$ is 
a \QI
invariant for such groups.  This yields a \QIc classification (\cite{ABC}, 
Corollary 5.6)
of these components of the \ABC groups.

In contrast, excluding endomorphisms with unipotent parts is essential to
the classification given here.  In fact, all nilpotent groups can be
expressed as the HNN extension of a simpler nilpotent group by a unipotent
matrix.  Thus, to classify \NBC groups defined with unipotent matrices
would be a very big step in the classification of nilpotent groups.  It
seems that we could use Pansu's invariants from \cite{P} as in \cite{ABC}
to any unipotent part which acts on the center $Z(G)$, but we will not
carry out that work here.

If an endomorphism fixes any vector, then it has unipotent part.  Thus,
any \NBC group which is an HNN extension defined by a unipotent-free
endomorphism has trivial center. However, a centerless \NBC group may be
defined by an endomorphism which is not unipotent-free if the unipotent
part acts on elements of $G$ which are not in the center of $G$.

\subsection{Permuted Absolute Jordan Form}\label{subsection:pajf}

The well-known Jordan form theorem states that any square matrix with
complex entries is conjugate to a matrix with the canonical Jordan form.  
In this subsection, we describe a modification of this form:
\emph{permuted absolute Jordan form}.

Let $\M_n(F)$ denote all $(n \times n)$-matrices over a field $F$, and let
$\GL_n(F)$ be the group of invertible matrices.

\begin{definition}
An matrix $J_n(\lambda) \in \M_n({\mathbb C})$ is a \emph{Jordan block}
with eigenvalue $\lambda$ if $J_n(\lambda)= (a_{ij})$ with
$$a_{ij} = \left\{ 
\begin{array}{cc}
\lambda & i=j\\
1       & i = j-1\\
0       & i \neq j-1,j.\\
\end{array}\right.$$
\end{definition}

That is, $J_n(\lambda)$ has $\lambda$'s along the diagonal, 1's on the
superdiagonal, and zeros elsewhere.

\begin{definition}
A matrix $J \in \M_n(\R)$ is a \emph{real Jordan block} if it has one of
the following two forms. The first form is an ordinary Jordan block
$J_n(\lambda)$ where $\lambda \in \R$. The second form, which requires
$n$ to be even, has a $2 \times 2$ block decomposition of the form
$$J = J_n(a,b) = 
\left(\begin{array}{ccccc}
Q(a,b) & \Id & \ldots & 0 & 0 \\
0 & Q(a,b) & \ldots & 0 & 0 \\  
\vdots & \vdots & \ddots & \vdots & \vdots \\
0 & 0 & \ldots & Q(a,b) & \Id \\
0 & 0 & \ldots & 0 & Q(a,b)
\end{array}\right)
$$
where $\Id$ is the identity, $0$ is the zero matrix,
$Q(a,b) = \textmatrix{a}{-b}{b}{a}$, and $b \ne 0$.
\end{definition}

\begin{definition} 
Given a partition $N = n_1 + \cdots + n_k$, and blocks
$J_i$ for $i = 1, \ldots, k$, we define the matrix 
$$M = (J_1, J_2, \ldots, J_k) = \left(\begin{array}{ccccc}
J_1 & 0 & \ldots & 0 \\
0 &  J_2 & \ldots & 0 \\  
\vdots & \vdots & \ddots &\vdots \\
0 & 0 & \ldots &  J_k\\
\end{array}\right).$$
If the blocks $J_i$ are all real Jordan blocks, then we say that matrix $M$ 
is in \emph{real Jordan block form}.  
If the blocks $J_i$ are all Jordan blocks, then we say that matrix $M$ 
is in \emph{complex Jordan block form}.  
If the blocks $J_i$ are all Jordan blocks, 
of the first type (that is, Jordan blocks with 
real $\lambda$), then we say that $M$ is in \emph{Jordan block form}.
\end{definition}

It is a standard result of linear algebra that every matrix $M \in
\M_n({\mathbb C})$ is conjugate via an element in $\GL_n({\mathbb C})$ to
a matrix in complex Jordan form, which is unique up to permutation of the 
Jordan blocks.  Similarly, every matrix $M \in \M_n(\R)$ is conjugate via an 
element in $\GL_n(\R)$ to a matrix in real Jordan form, which is also unique 
up to permutation of the Jordan blocks.

\begin{definition}
The \emph{absolute Jordan form} of $M \in \M_n(\R)$ is the matrix obtained
by replacing each diagonal entry of the complex Jordan form with its 
absolute value.
\end{definition}

We will resolve the nonuniqueness in absolute Jordan form by specifying: 
\begin{itemize}
\item if $|\lambda_i|>|\lambda_j|$, then $i > j$, and 
\item if $|\lambda_i|=|\lambda_j|$, and $n_i > n_j$, then $i > j$. 
\end{itemize}

\begin{definition}
Suppose $<$ defines a partial order on the set $\N_n = {1, \ldots, n}$ 
We say a permutation $\sigma \in S_n$ \emph{preserves the partial order 
$<$} if $$a < b \implies \sigma(a) < \sigma(b).$$
\end{definition}

\begin{definition}
If a matrix $M = (J_{n_1}(\lambda_1), J_{n_2}(\lambda_2), \ldots,
J_{n_k}(\lambda_k))$ is in Jordan form, let $N_0 = 0$ and $N_i = n_1 + 
\cdots + n_{i-1}$.  
Define a partial order $<_M$ on $\N_n$ by: $i <_M j$ if:
\begin{itemize} 
\item $i < j$  as natural numbers, and 
\item $N_{l-1} < i,j \leq N_{l}$ for some $l = 1, \ldots, k$.
\end{itemize}
If there is no $l$ such that $N_{l-1} < i,j \leq N_{l}$, then $i$ and $j$
are not comparable by $<_M$. We say $<_M$ is the \emph{partition
associated to the Jordan matrix $M$}. \end{definition}

\begin{definition} A matrix $M = (m_{i,j})$ is in \emph{permuted Jordan
form} if there is a permutation $\sigma \in S_n$ such that 
\begin{itemize}
\item the matrix
$M_\sigma = (m_{\sigma(i),\sigma(j)})$ is in Jordan form, and 
\item the permutation $\sigma^{-1}$ preserves the partial order $<_M$ 
associated to $M_\sigma$.
\end{itemize}
\end{definition}

The diagonal entries of $M_\sigma$ will be the diagonal elements of $M$,
permuted by $\sigma$.  As a consequence of the second condition of the
definition, all the off-diagonal 1's will be above the diagonal of matrix
$M$.

Given an endomorphism $\phi$ of a Carnot group $G$, we associate to this a 
canonical matrix in permuted absolute Jordan form as follows:

\begin{definition} Consider an endomorphism $\phi$ of an $n$-dimensional
nilpotent Lie group $G$.  The induced linear map of the Lie algebra can be
represented with respect to some basis $\{e_i\}$ by a matrix $M$ in Jordan
form. By Proposition \ref{theorem:linalg} this basis respects the 
nilpotent
grading.  Therefore, Theorem \ref{theorem:kar} associates a set of integer
weights $\{w_i\}$ to the basis $\{e_i\}$. We define a permutation $\sigma 
\in S_n$ so that: 
\begin{itemize} 
\item if $w_i > w_j$, then $\sigma(i) < \sigma(j)$, and
\item if $w_i=w_j$, then $\sigma(i)<\sigma(j)$ if and only if $i<j$. 
\end{itemize} 
\end{definition}

The effect of these requirements is that we do as little rearranging of 
the Jordan form as possible, subject to the condition that the new 
weight vector is non-increasing.

%\example

%COPIED FROM DRAFT 423
%\example $$M = (J_3(2), J_3(3)) = \left[ \begin{array}{cccccc} 
%2&0&0&0&0&0\\
%1&2&0&0&0&0\\ 0&1&2&0&0&0\\ 0&0&0&3&0&0\\ 0&0&0&1&3&0\\ 0&0&0&0&1&3\\
%\end{array} \right] $$
%with weights:  $w_1 = 1$, $w_2 = 2$, $w_3 = 3$, $w_4 = 1$, $w_5 = 1$, $w_6
%= 2$.

%$$[f_1] = \max\{2^t, (t \cdot 2^t)^{\frac{1}{2}}, (t^2 \cdot
%2^t)^{\frac{1}{3}}\} = 2^t$$$$
%[f_2] = \max\{(2^t)^{\frac{1}{2}}, (t \cdot 2^t)^{\frac{1}{3}}\}
%=2^{\frac{t}{2}}$$$$
%[f_3] = \max\{(2^t)^{\frac{1}{3}}\}=2^{\frac{t}{3}}$$$$
%[f_4] = \max\{3^t, t \cdot 3^t, (t^2 \cdot 3^t)^{\frac{1}{2}} \} =
%t \cdot 3^t$$$$
%[f_5] =  \max\{3^t, (t \cdot 3^t)^{\frac{1}{2}}\}= 3^t$$$$
%[f_6] =\max\{(3^t)^{\frac{1}{2}}\} =3^{\frac{t}{2}}$$

%Thus, the set of divergence rates $div(\Gamma) = \{t \cdot 3^t, 3^t, 2^t,
%3^{\frac{t}{2}}, 2^{\frac{t}{2}}, 2^{\frac{t}{3}}\}$. This is an
%interesting example, because there are six
%distinguishable growth rates, as many as there are dimensions of the
%space.  

Define a three-step nilpotent group 

Let $N = H \times H \times H$, where $H$ is the Heisenberg group.  Thus, 
$N$ is a two-step nilpotent group with presentation:
$$N = \langle a_1, a_2, a_3, a_4, a_5, a_6, a_7, a_8, a_9 \suchthat [a_1, a_2] 
= a_3, [a_4, a_5] = a_6, [a_7, a_8] = a_9 \rangle.$$

Consider the endomorphism of the associated Lie group 
$\g$ which is represented by the matrix:

$$M = \left( \begin{array}{cccccc} 
2&1&0&0&0&0\\
0&2&0&0&0&0\\ 
0&0&2&0&0&0\\ 
0&0&0&3&1&0\\ 
0&0&0&0&3&1\\ 
0&0&0&0&0&3\\
\end{array} \right),$$
with associated weight vector $w = (1, 2, 3, 1, 1, 2)$.  

Then the 
permuted absolute Jordan form will be obtained with the permutation 
$$\sigma = \left( \begin{array}{cccccc}  
1&2&3&4&5&6\\ 
1&4&6&2&3&5\\
\end{array} \right).$$

The permuted absolute Jordan matrix will be:
$$M_\sigma = \left(\begin{array}{cccccc} 
2&0&0&1&0&0\\
0&3&1&0&0&0\\ 
0&0&3&0&1&0\\ 
0&0&0&2&0&1\\ 
0&0&0&0&3&0\\ 
0&0&0&0&0&2\\
\end{array} \right).$$
\bigskip

\subsection{Linear Algebra}\label{subsection:linalg}

In order to understand how the induced map $\phi^*$ acts on $\g$, we need 
the following results of linear algebra.  Recall that $\phi$ weakly 
preserves a decomposition $V = L_1 \oplus \cdots \oplus L_k$ 
if $\phi( L_i \oplus \cdots \oplus L_k) \subset L_i \oplus \cdots \oplus 
L_k$ for all $i = 1, \ldots, k$.

\begin{prop}\label{theorem:linalg}
Let $\phi$ be a nonsingular linear map of a vector space $V$.  
Suppose that $V$ has a direct sum decomposition
$V = L_1 \oplus \cdots \oplus L_k$ 
which is weakly preserved by $\phi$.
Then there is a basis $\{e_i\}$ for $V$ such that 
\begin{enumerate}
\renewcommand{\labelenumi}{\alph{enumi})}
\item  The matrix representing $\phi$ in the basis $\{e_1, \ldots, e_n\}$
is in permuted Jordan form.
\item The permutation which puts the matrix in permuted Jordan form 
respects the partial order associated to the Jordan form.  (Equivalently, 
the matrix is in upper triangular form)
\item If $d(i) = \dim (L_i \oplus \cdots \oplus L_k)$, then 
$L_i \oplus \cdots \oplus L_k = \spn \{e_1, \ldots, 
e_{d(i)}\}$
\end{enumerate}
\end{prop}

The proof of the theorem will rely upon the following lemma, in which we 
consider the special case that $\phi$ acts as a single Jordan block.

\begin{lemma}\label{lemma:linalg}
Suppose $\phi$ is a linear map of an $n$-dimensional vector space $V$, and 
that $\phi$ can be represented as a matrix which is a single 
real Jordan block with respect to the 
basis $\{e_1, \ldots, e_n\}$.  
Let $V_i = \spn \{e_1, \ldots, e_i\}$.
If $W \subset V$ is a subspace preserved by $\phi$, (i.e., $\phi(W) 
\subset W$), then $W=V_{\dim(W)}$. 
\end{lemma}

\begin{proof} 

\textbf{Case 1:} 
The real Jordan block representing $\phi$ is 
$J_n(\lambda)$.  

Suppose $W\neq V_i$ for all $i$.  Let $i$ be maximized subject to the
condition that $V_i \subset W$.  Yet, by assumption, $W \neq V_i$, so we can
choose some $w \in W \cap \spn\{e_{i+1}, \ldots, e_n\}$, which can be
expressed as $w = a_{i+1}e_{i+1} + \cdots + a_ne_n$.  Let $j$ be maximized
subject to $a_j \neq 0$; that is, $w \in W \cap \spn\{e_{i+1}, \ldots,
e_j\}$.  We will demonstrate that there is some $w' \in W \cap
\spn\{e_{i+1}, \ldots, e_{j-1}\}$.

This will imply, by induction, that there is some element $w^* \in W \cap 
\spn\{e_{i+1}\}$,
which will imply that $V_{i+1} \subset W$. This contradiction will imply
that $W = V_i$, for some $i$.  Then, $i$ is determined so that
$\dim(W)=\dim(V_i)=i$.

To construct $w'$, note that $\phi(w) = \lambda w +
a_{2}e_{1} + \cdots + a_{i}e_{i-1} + a_{i+1}e_{i} + \cdots + a_j e_{j-1}$.  
But $\lambda w \in W$ and $a_{2}e_{1} + \cdots + a_{i}e_{i-1} 
+ a_{i+1}e_i \in V_i \subset W$, so $w' = 
a_{i+2}e_{i+1} + \cdots +
a_je_{j-1} \in W$ as well.

\textbf{Case 2:} 
The real Jordan block representing $\phi$ is 
$J_n(a,b)$ (and therefore $n = 2m$ is even).

This case proceeds similarly to the first, except that elements in the
vectors are replaced by ordered pairs.  Given an element $w \in \R^n$,
express $w = ((c_1, d_1), \ldots, (c_m,d_m))$, where each $(c_i,d_i) \in
\R^2$.

Again, assume that $W \neq V_i$ for all $i$.  Let $i$ be maximized subject
to the condition that $V_{2i} \subset W$.  Yet, by assumption, $W \neq
V_{2i}$, so we can choose some $w \in W \cap \spn\{e_{i+1}, \ldots, 
e_n\}$,
which can be expressed as $w = c_{i+1}e_{2i+1} +d_{i+1}e_{2i+2} + \cdots + 
c_{m}e_{n-1} +d_{m}e_{n}$.  Let $j$
be maximized subject to $(c_j,d_j) \neq (0,0)$; that is, $w \in W \cap 
\spn\{e_{2i+1}, \ldots, e_{2j}\}$.  We will demonstrate that there is 
some 
$w' \in W \cap \spn\{e_{2i+1}, \ldots, e_{2j-2}\}$.

But, by induction, this implies that there is some element $w^* \in W \cap 
\spn\{e_{i+1}, e_{i+2}\}$.  Furthermore, $\phi(w^*)  \in W \cap 
\spn\{e_{i-1}, e_i, e_{i+1}, e_{i+2}\}$.  Since   
$\spn\{e_{i-1}, e_i\} \subset W$, we can subtract the $e_{i-1}$ and 
$e_i$ components from $\phi(w^*)$.  The result is an element of 
$\spn\{e_{i+1}, e_{i+2}\}$ which is linearly independent of 
Since $w^*$.  This implies that $V_{i+2} \subset W$. This contradiction 
will imply that $W = V_i$, for some
$i$.  Then, $i$ is determined so that $\dim(W)=\dim(V_i)=i$.

Now we construct $w'$. Let $Q = Q(a,b)$, so that  
$$\phi(w) = (Q(c_1,d_1) + (c_2, d_2), Q(c_2, d_2) + (c_3, d_3), \ldots, 
Q(c_m, d_m)).$$ 
Since $(c_m, d_m)$ and $Q(c_m, d_m)$ are linearly independent vectors in 
$\R^2$, there are constants $r, s \in \R$ such that 
$r (c_m, d_m) + s Q(c_m, d_m) = 0$.  Thus, 
$rw + s\phi(w) \in W \cap \spn\{e_1, \ldots, e_{2j-2}\}$.  Since  $V_{2i} 
\subset W$, we can subtract off the neccesary terms of $rw + s\phi(w)$ to 
obtain $w^*  \in W \cap \spn\{e_{i+1}, \ldots, e_{2j-2}\}$.
\end{proof}

\begin{proof}[Proof of Proposition \ref{theorem:linalg}]
It is a standard result of linear algebra (Jordan Decomposition Theorem) 
that every nonsingular linear map on a vector space can be decomposed 
into the action on a direct sum of $\phi$-invariant root spaces, such 
that the induced action on each root space is as a single real Jordan block.  

Let the root space decomposition of $V$ induced by $\phi$ be $V = R_1
\oplus \cdots \oplus R_m$.  On each $R_i$, the map $\phi$ can be
represented as a single real Jordan block with respect to some basis for 
$R_i$, say,
$\{a^i_1, \ldots, a^i_{\dim(R_i)}\}$.  Let
$V^i_j = R_i \cap (L_j \oplus \cdots \oplus L_k)$.  Since $\phi$
preserves both $R_i$ and $L_j \oplus \cdots \oplus L_k$, it also
preserves $V^i_j$.  Thus, we apply Lemma 9 to prove that 
$$V^i_j  = \spn \{a_1^i, \ldots, a^i_{\dim(V^i_j)}\}.$$ 

Define a function $f\colon \{a_j^i\} \rightarrow\{1, \ldots, n\}$ by 
listing the $a_j^i$ as follows:
$$a^1_1, \ldots,  a^1_{\dim V_1^1}, a_1^2, \ldots, a^2_{\dim V_1^2}, 
\ldots, 
a_1^m, \ldots, a^m_{\dim V_1^m},\ldots $$
$$a^i_{\dim V_1^1 +1}, \ldots,  a^i_{\dim V_2^1}, \ldots, a^m_{\dim 
V_k^m}.$$
Then, $f$ maps the $i^{th}$ element of this list to $i$.
Thus, the first $\dim V^1_1 + \cdots + \dim V_1^m = \dim L_1$ vectors 
span 
$L_1$, and similarly for $L_j$.  Also, for fixed $i$, the order of 
$\{a_j^i\}$ is preserved.  Thus, the matrix representing $\phi$ with 
respect to the basis $\{e_i\}$
is in permuted Jordan form.  
Furthermore, this permutation respects the partial order associated to 
the Jordan matrix.
\end{proof}

\subsection{An Order on Divergence Rates} \label{subsection:order}
In order to compare the rates at which the lengths of different vectors 
grow, we will describe a partial order and an equivalence relation on
functions from $\R$ to $\R$.  The order will characterize functions 
which dominate others.

\begin{definition} 
Given $f, g \colon \R \rightarrow \R$ we say that $f
\preceq g$ if there are $K, C \geq 0$ such that

$$f(t) \leq Kg(t) + C \mbox{ for all } t \in \R .$$ 

We say $f, g$ are \emph{comparable} ($f \simeq g$) if $f \preceq g$ and $g 
\preceq f$.  
\end{definition}

Equivalently, $f\simeq g$ if there exist constants $K', 
C' > 0$ such that:
$$\frac{1}{K'} g(t) -C'  \leq f(t) \leq K'g(t) +C' \mbox{ for all } t \in 
\R.$$

This is an equivalence relation: 
\begin{itemize}
\item {(Symmetry)} 
For all $f$, $f \simeq f$ with constants $K=1$, $C=0$ 
\item {(Transitivity)} If $f \simeq g$ with constants $K$ and $C$, and $g 
\simeq h$ with constants $K'$ and 
$C'$, then $f \simeq h$ for constants $KK'$ and $\max \{K'C +C, 
\frac{C}{K'} + C'\}$
\item {(Reflexivity)} If $f \simeq g$ with constants $K$ and $C$, then $g 
\simeq f$ with constants $K$ and $\max \{CK, \frac{C}{K}\}$.
\end{itemize}

We denote the equivalence class of a function $f$ by
$[f]$. In the \ABC case, the corresponding divergence rates are always
exponential or polynomial$\cdot$exponential functions.  In Theorem 
\ref{theorem:divrates}, we will show that divergence rates are
exponential, polynomial$\cdot$exponential, or roots thereof.  For this
class of functions, the partial order is actually an order.

\begin{lemma}
\label{lemma:class}
For the class of functions $\C = \{
f(t) = (t^n \cdot \lambda^t)^{\tfrac{1}{d}} 
| \lambda \neq 1, \lambda \geq 0, n,d, \in \Z^+\}$, the partial 
order $\preceq$ is an order.
\end{lemma}

\begin{proof}
Consider two such functions
$f(t) = (t^n \cdot \lambda^t)^{\frac{1}{d}} $
and $g(t) = (t^m \cdot {\kappa}^t)^{\frac{1}{c}}$ with $\lambda, \kappa > 
0, n, m, c, d \in \Z^+$.

\textbf{Case 1:} 
$\lambda^{\frac{1}{d}} \neq {\kappa}^{\frac{1}{c}}$

In this case, the exponential growth dominates.  Assume, without loss of generality,
that $\lambda^{\frac{1}{d}} > {\kappa}^{\frac{1}{c}}$.  If we set $C=0$, 
we are looking for a value of $K$ such that 
$$\frac{(t^m \cdot \kappa^t)^{\frac{1}{c}}}{(t^n \cdot \lambda^t)^{\frac{1}{d}}} 
\leq K \mbox{ for all } t.$$
Elementary calculus shows
that $f(t) = t^ab^t$ has a global maximum at $t=\frac{-a}{\ln(b)}$, where 
$a=\frac{m}{c} - \frac{n}{d}$
and $b=\frac{{\kappa}^{\frac{1}{c}}}{\lambda^{\frac{1}{d}}} $.  Thus, choosing
for $K$ the maximal value of $f$ shows that $g \preceq f$.

\textbf{Case 2:} 
$\lambda^{\frac{1}{d}} = {\kappa}^{\frac{1}{c}}$

In this case, the exponential parts of the functions grow at the same 
rate, so the degree of polynomial determines which function is larger.  
Assume, without loss of generality, that $\frac{n}{d} > \frac{m}{c}$.
For all $t \geq 1$, $f(t) \geq g(t)$.  
Since $g(t) < \lambda$ for all $0 < t < 1$, we choose $C=\lambda$, $K=1$, 
to show that $g \preceq f$.
\end{proof}

We will also need the following Lemma.

\begin{Lemma}
\label{lemma:sum}
Suppose $\{f_1, \ldots, f_n\} \subset \C$ has $f_1 \succeq f_i$ for $i = 
2, \ldots, n$.  Then 
$[\Sigma_{i=1}^n f_i] = [f_1]$.
\end{Lemma}

\begin{proof} 
If for each $i=2, \ldots, n$, $f_1 \succeq f_i$ with constants $K_i, C_i$, then $[\Sigma_{i=1}^n f_i]
= [f_1]$ with constants $1+ K_2 + \cdots + K_n$ and $C_2 + \cdots + C_n$.
\end{proof}

\subsection{Growth Rates of Vectors}
\label{subsection:vectorgrowth}
Given a 1-parameter subgroup $M^t$ of $\GL_n(\R)$, Farb and Mosher
derive \cite{ABC} upper and lower bounds for the growth of vectors 
$||M^tv||$.  We
will need more: bounds on the growth of each coordinate of $M^tv$.  In the
special case of $v=e_i$, one of the standard basis vectors for $\R^n$,
results are actually contained in their proof and are stated below in
Proposition \ref{theorem:vectorgrowth1}. In Theorem
\ref{theorem:vectorgrowth2} we extend to a result for arbitrary $v \in 
\R^n$.  We will not give the growth function explicitly, but only up to 
the relation $\simeq$ defined in Subsection \ref{subsection:order}.  

\begin{prop}\label{theorem:vectorgrowth1}
Consider the Jordan block
$M= J_{n}(\lambda) \in \GL_n(\R)$
(In particular, $\lambda \in \R$.)  
Let $\{e_i\}$ be the standard basis for $\R^N$.  
Then, considered as functions of $t$:
$$e_k \cdot M^t \cdot e_j \simeq 
\left\{ 
\begin{array}{cc}
\lambda^t t^{k-j}  &   j \leq k \leq n\\
0 & k < j.
\end{array}\right.$$
\end{prop}

\begin{proof}
See \cite{ABC}, Equation 3.1.
\end{proof}

Now we extend this result to an arbitrary $x \in \R^n$:

\begin{prop}
\label{theorem:vectorgrowth2}
Consider $x = (x_1, \ldots, x_n) \in \R^n$ and $M \in \GL_n(\R)$ a Jordan 
matrix.  As functions of 
$t$:
$$e_k \cdot M^t \cdot x \simeq \max_i 	
\{e_k \cdot M^t \cdot e_i| x_i \neq 0\}.$$
\end{prop}

Although the equivalence class of $e_k \cdot M^t \cdot x$ is 
independent of the nonzero values of $x_i$, the constants implicit in 
the relation $\simeq$ does depend linearly on $x_i$ for some $i$.

\begin{proof} We write $x \in \R^n$ as $ x = \Sigma_{i=1}^n x_ie_i$.  
Then, $$M^t \cdot x = M^t \cdot \Sigma_{i=1}^n x_i e_i = \Sigma_{i=1}^n
x_i M^t e_i.$$ Restricting our attention to the $k^{th}$ coordinate: $$e_k
\cdot M^t \cdot x = \Sigma_{i=1}^n x_i e_k M^t e_i.$$ Now, we consider the
divergence rate (as a function of $t$) of this coordinate, using the
equivalence relation described in Subsection \ref{subsection:order}.  
Proposition
\ref{theorem:vectorgrowth1} tells us that each of the component functions
$e_k M^t e_i$ is in $\C$, and so, by Lemma \ref{lemma:class}, all are
comparable.  Lemma \ref{lemma:sum} implies that the sum of such functions
is equivalent to the maximum.  For fixed $x_i \neq 0$, we have $x_i f(t) 
\simeq
f(t)$ Thus, $$e_k \cdot M^t \cdot x \simeq \max_{i=1, \ldots, n} \{e_k
\cdot M^t \cdot e_i \suchthat x_i \neq 0\}.$$ \end{proof}

\subsection{Putting the Pieces Together: Divergence Rates} 
\label{subsection:divrates} Given a
connected, simply-connected nilpotent Lie group $G$ and an injective
endomorphism $\phi$ of $G$, we get an induced linear map $\phi^*$ on the
Lie algebra $\g$.  Assume this can be represented by a Jordan matrix
$M = (J_{n_1}(\lambda_1), J_{n_2}(\lambda_2), \ldots,
J_{n_k}(\lambda_k))$ with respect to some basis $\{e_i\}$.  
Let $N_i = n_1 + \cdots + n_i$.
By Lemma \ref{theorem:endo},
this basis is also consistent with the nilpotent grading of $\g$ (although
perhaps in a permuted order).

We define a Lie group $G_\phi = G \rtimes_\phi \R$, where the action of
$\R$ is given by $\phi^t$.  We identify $G_\phi$ with $\R^n \times \R$.  
Given any left-invariant metric on $G$, the metric on each slice $G_t = G
\times \{t\}$ is given by the pullback of $\phi^t \colon G \rightarrow G$.
See Subsection \ref{subsection:model} for a description of the group operation.

In this subsection, we will calculate the rate of divergence of ``vertical
flow lines'' in $G_\phi$; that is, for a given $x \in \R^n$, we calculate
the distance from $(0,t)$ to $(x,t)$ within $G_t$, which we denote by 
$f_x(t)$.

Note that a multiset is similar to a set in that the elements do not
have a designated order; it differs from a set in that elements may be
repeated.

\begin{definition}
Given a Lie group $G_\phi$ as above, and a basis $\{e_i\}$ 
so that $\phi^*$ is represented as a real Jordan matrix, let
$f_i(t) = d_t((0,t),(e_i,t))$.  Recall that $[f]$ denotes
the equivalence class of $f \colon\R \rightarrow \R$ as 
defined in Subsection \ref{subsection:order}.  We define a 
multiset of divergence rates:
$$\D_\phi = \{[f_i(t)]\}.$$
\end{definition}

\begin{theorem}
%[Divergence Rates are Elements of $\C$]
\label{theorem:divrates}
The multiset of divergence rates for a Lie group $G_\phi$ is:
\begin{enumerate}
\renewcommand{\labelenumi}{\alph{enumi})}
\item independent of generating set $\{e_i\}$;
\item a finite subset of $\C$ (as defined in Lemma \ref{lemma:sum}) ; and
\item contains at least one element of the form $[\lambda^t]$.
\end{enumerate}
\end{theorem}

The proof of this theorem will require two lemmas.
As in the \ABC case, we will see that $f_x(t) = || M^{-t}(x) ||$, although
this requires more work.  On the other hand, now the norm $|| \cdot ||$ is
the norm in the nilpotent metric, which depends on the grading of the
nilpotent group, as discussed in Subsection \ref{subsection:kar}: 
$$||(x_1, x_2, \ldots, x_n)|| \sim \max_{i} \{ x_i^{\frac{1}{w_i}}\},$$
where $w_i =$ weight of $x_i$ in the grading of $N$. 

\begin{lemma} 
\label{lemma:divrates1}
With $G_\phi$, $\{e_i\}$, and $f_i(t)$ given as above, reorder the 
$\{e_i\}$ so that 
$$f_1(t) \preceq f_2(t) \preceq \cdots \preceq f_n(t),$$
and define $V_i = \spn \{e_1, \ldots, e_i\}$.  If $x \in V_i 
\setminus V_{i-1}$, then $[f_x(t)]= [f_i(t)]$.
\end{lemma}

Notice that, although we write $x$ as a vector, it is not a vector in the
Lie algebra, but, rather, a fixed point in $\R^n$, which 
corresponds to a (varying) point in the Lie group, expressed in
the coordinates corresponding to our chosen basis for the Lie algebra.

\begin{proof}
Let $\{e_i\}$ be an orthonormal basis for the metric on $G_0 = G \times
\{0\}$.  Then the metric on $G_t$ can be defined equivalently by the
orthonormal basis $e_i(t) = (\phi^t)^*(e_i) = M^t \cdot e_i$. The results
described in Subsection \ref{subsection:kar} apply to each slice $G_t$ to show 
that balls in this left-invariant metric are comparable to polynomial
ellipsoids.  The proof in \cite{K} shows more. Because the commutivity
data for $\{e_i(t)\}$ is independent of $t$, the weights $w_i$ and
constant $a$ in Theorem \ref{theorem:kar} are the same for each slice
$G_t$.  Thus, it suffices to write $x$ in terms of this basis; i.e., $x =
\Sigma_i x_i(t) e_i(t)$.  This is simply a change of basis, so the
coordinates $(x_1(t), \ldots, x_n(t))$ are given by the vector $M^{-t}
\cdot x$.   That is, $$d_t((0,t),(x,t)) = ||M^{-t}x||.$$

By Corollary \ref{cor:dist1},
$$||M^{-t}x|| = \max_k \{|e_k \cdot M^{-t}\cdot x|^{\frac{1}{w_k}} \}.$$
Applying Proposition \ref{theorem:vectorgrowth2},
\begin{equation}
\begin{split}
||M^{-t}x|| &= \max_k \max_i \{|e_k \cdot M^{-t}\cdot 
e_i|^{\frac{1}{w_k}}\suchthat x_i \neq 0 \}\\
&=\max_i \{ f_i(t) \suchthat x_i \neq 0 \}\\
&=f_{\max\{i \suchthat x_i \neq 0\}}(t).\\
\end{split}
\end{equation}
Thus, for $x \in V_i \setminus V_{i-1}$, we have shown that $||M^{-t}x|| 
\simeq f_i(t)$.
\end{proof}

\begin{proof}[Proof of Theorem \ref{theorem:divrates}(a)]
Suppose that for a different Jordan basis $\{e_i'\}$, we obtain a different multiset of divergence rates: 
$$g_1(t) \preceq g_2(t) \preceq \cdots \preceq g_n(t).$$
Let $i$ be maximized subject to $g_i \not\simeq f_i$ and assume,
without loss of generality, that $g_i \prec f_i$.  Let $V_i = \{x \in G
\suchthat f_x \preceq g_i\}$.  Then, $V_i \supset \spn \{e_1', \ldots,
e_i'\}$ and so $\dim(V_i) \geq i$.  On the other hand, $\spn \{e_i,
\ldots, e_n\} \cap V_i = \emptyset$, so $\dim V_i < i$.  This
contradiction shows that $f_i = g_i$ for all $i$.
\end{proof}

\begin{lemma}
\label{lemma:divrates2}
With the notation as in Proposition \ref{theorem:vectorgrowth1}, 
let $f_j(t) = f_{e_j}(t)$, and suppose $e_k$ has weight $w_k$ given by 
Theorem \ref{theorem:kar}.  Then,
$$f_j(t) \simeq \max_k \{(t^{(k-j)}\lambda_i^t)^{\frac{1}{w_k}}| j 
\leq k \leq N_i \}.$$
\end{lemma}

\begin{proof}
We simply apply the metric determined in Theorem \ref{theorem:kar} to the 
coordinates determined by Proposition \ref{theorem:vectorgrowth1}.
\end{proof}

\begin{proof}[Proof of Theorem \ref{theorem:divrates}(b) and (c)]
Lemma \ref{lemma:divrates2} implies that each divergence rate $f_i \in \C$, which is (b).  

For $j = N_i$ (corresponding to $e_j$ is an eigenvector) we have $f_j(t) = 
(\lambda_j^{\frac{1}{w_j}})^t$, which establishes (c).
\end{proof}

\example Consider the endomorphism $\phi$ of the Heisenberg group $H$
defined by \begin{equation*} \begin{split} \phi(x) &=x^3yz\\ \phi(y)
&=x^{-1}y\\ 
\end{split} \end{equation*} 
Together these imply:$ \phi(z)=z^4$, and the endomorphism is represented in this basis by the matrix
$$M = \left( \begin{array}{ccc} 3&-1&0\\ 1&1&0\\ 1&0&4\\ \end{array}
\right),$$ with associated weight vector $w = (1, 1, 2)$.  The permuted
absolute Jordan form of this matrix is $$M' = \left( \begin{array}{ccc}
3&0&0\\ 0&2&1\\ 0&0&2\\ \end{array} \right).$$

The associated divergence rates are then:
\begin{equation*}
\begin{split}
[f_z] &= (3^t)^{\frac{1}{2}}\\
[f_x] &= t \cdot 2^t\\
[f_y] &= 2^t
\end{split} \end{equation*} 

%Thus, the set of divergence rates $div(\Gamma) = \{t \cdot 3^t, 3^t, 2^t,
%3^{\frac{t}{2}}, 2^{\frac{t}{2}}, 2^{\frac{t}{3}}\}$. This is an
%interesting example, because there are six
%distinguishable growth rates, as many as there are dimensions of the
%space.  
\bigskip

\section{Proof of the Classification}
\label{section:coarsetop}
%Section 3
This section contains the proofs of three of the main results of this 
paper (Theorems \ref{theorem:pajf}-\ref{theorem:invariant2}), which 
together constitute significant progress towards a 
classification of nonpolycyclic \NBC groups.  Subsection \ref{subsection:model} 
describes a geometric model space for the groups.  Subsections 
\ref{subsection:opjs} and \ref{subsection:reducetoajf} contain the proof of 
Theorem \ref{theorem:pajf}.  Subsections 
\ref{subsection:coarsetop}-\ref{subsection:time2} 
contains the proof of Theorem \ref{theorem:invariant1}.  The proof of 
Theorem \ref{theorem:invariant2} is completed in Subsection 
\ref{subsection:growth}.

\subsection{Model Spaces for Nilpotent-by-Cyclic Groups}
\label{subsection:model}

Given a discrete nilpotent group $N$ which is a lattice in a nilpotent Lie
group $G$ and an injective endomorphism $\phi$ of $N$ (which always
extends to an endomorphism of $G$), recall that $\Gamma_\phi$ is the HNN
extension of $N$ by $\phi$, and the Lie group $G_{\phi}$ is the
semi-direct product $G \rtimes_\phi \R = \{(x,t) \suchthat x \in G, t \in
\R\}$.  Multiplication is defined by:  $$(x,t) \cdot (y,s) = (x \cdot_G
\phi^t(y) , t+s)$$ for all $(x,t), (y,s) \in G \times \R$.  Recall that
the left-invariant metric on $G_\phi$ is defined in Subsection
\ref{subsection:divrates}.

In this subsection, we construct a metric complex $X_{\phi}$ on which 
$\Gamma_{\phi}$
acts properly discontinuously and cocompactly by isometries.  Thus,
$\Gamma_{\phi}$ will be \QIc to the metric space $X_{\phi}$.  This
construction will parallel the presentation in \cite{ABC}.

Let $M$ be the $n$-manifold with fundamental group $\pi_1(M) =N$.
Since $\Gamma_{\phi}$ is an ascending HNN extension of $N$, 
it is the fundamental group of the mapping torus of $M$ under the 
endomorphism $\phi$.  Let $X_{\phi}$ be the universal cover of this 
mapping torus.  Topologically, $X_{\phi} \approx \R^n \times T_\phi$, 
where $T = T_\phi$ is the Bass-Serre tree associated to the HNN extension.  
Thus, $T$ is a homogenous directed tree with one edge oriented in and 
$[N:\phi(N)]$ edges oriented out of each vertex.

Then $X_{\phi}$ is a fiber product of $G_{\phi}$ and $T_\phi$ over
$\R$.  The Lie group $G_\phi$ comes naturally equipped with a height
function $h((x,t)) = t$.  We define a height function on $T_\phi$ as
follows. Fix a base point $x$ in the tree $T_\phi$, and define a path
metric on the tree such that each edge has unit length.  This gives a
height function $h \colon T_\phi \rightarrow \R$, defined by $|h(y)| =
d(x,y)$ and $h(y) > 0$ if and only if the distance minimizing path from
$x$ to $y$ begins with an edge which is oriented out of $x$.  
The metric on $X_\phi$ is defined by the fiber product of the 
metrics on $T_{\phi}$ and $G_{\phi}$.

There are induced projections $g_\phi \colon X_{\phi} \to G_{\phi}$ and
$\pi_\phi \colon X_\phi \to T_\phi$, and an induced height function
$X_\phi \to \R$.

A \emph{horizontal leaf} $L \subset X_{\phi}$ is a subset of the form $L =
\pi_\phi^{-1}(x)$ where $x \in T_\phi$.  Let $\ell$ be a bi-infinite line
in the tree $T_\phi$.  Then, a \emph{hyperplane} $P_\ell \subset X_\phi$ 
is a subset of the form $P_\ell = \pi_\phi^{-1}(\ell)$.
If the line $\ell$ is coherently oriented in $T_\phi$, then $P_\ell$ is
isometric to $G_\phi$ by construction, and we call $P_\ell$ a {\em 
coherent hyperplane} in
$X_\phi$.  If the line $\ell$ is not coherently oriented in $T_\phi$ (and
thus switches orientation exactly once), then we call $P_\ell$ a {\em 
incoherent hyperplane} in $X_\phi$.  We will show, in Theorem 
\ref{theorem:bigon}, that such a hyperplane is not \QIc to $G_\phi$.

Suppose that $\F$ is a decomposition of a metric space $X$ into disjoint
subsets whose union is $X$.  Let $\G$ be such a decomposition of a metric
space $Y$.  A quasi-isometry $f \colon X \rightarrow Y$ \emph{coarsely
respects} the decompositions $\F$ and $\G$ if there exists an $A \geq 0$
and a map $h \colon \F \rightarrow \G$ such that for each element $L \in
\F$ we have $d_\Haus(f(L), h(L)) \leq A$.  
For example, we will refer to quasi-isometries which \emph{coarsely 
respect horizontal leaves}  or \emph{coarsely respect 
vertical flow lines}.

\subsection{One Parameter Jordan Subgroups} 
\label{subsection:opjs}

In this subsection and the next, we will prove Theorem \ref{theorem:pajf}.  
Consider the geometric model spaces associated to two \NBC groups with the
same permuted absolute Jordan form.  Since these model spaces are 
fiber products, it will suffice to show that (1) the
associated Lie groups are \qic, (2) the associated trees (of which the 
model spaces are fiber products) are \qic, and (3) the quasi-isometries 
between them have induced time change functions which are consistent.
In this subsection, we will set the stage for step (1) by establishing a
relationship between a matrix and its permuted absolute Jordan form.

Given a matrix $M\in \GL_n(\R)$ in Jordan form (not just in real Jordan
form--no $J_n(a,b)$ blocks), we say that $\rho(t) = e^{Mt}$ is a
1-parameter \emph{Jordan subgroup}.  The matrices $e^{Mt}$ may not be in
Jordan form. Not only can a single matrix be conjugated into Jordan form,
but Witte has shown (\cite{ABC}, Theorem 3.1, also \cite{Wit}) that an
entire 1-parameter subgroup of $\GL_n(\R)$ can be transformed into a
1-parameter Jordan subgroup.  We will need a corollary of this theorem,
under the additional hypothesis that $M$ weakly preserves a grading of
$\R^n$.

\begin{theorem}[\cite{ABC}, Theorem 3.1, 1-parameter real Jordan form]
\label{theorem:ajf} Let $M^t$ be a 1-parameter subgroup of $\GL_n(\R)$.  
There exists a 1-parameter Jordan subgroup $e^{Jt}$, a matrix $A \in
\GL_n( \R)$ and a bounded 1-parameter subgroup $P^t$ conjugate into the
orthogonal group O$(n,\R)$, such that $e^J$ is the absolute Jordan form of
$M$, and letting $\barM^t = A^{-1}e^{Jt}A$ we have $$M^t = \barM^t P^t =
P^t\barM^t.$$ \end{theorem}

\begin{proof}(from \cite{ABC}, p. 156) Given a
general 1-parameter subgroup $e^{\mu t}$ in $\GL_n(\R)$, choose $A$ so that
$A^{-1} \mu A$ is in real Jordan form, and so $A^{-1} \mu A = \delta + \nu +
\eta$ where $\delta$ is diagonal, $\nu$ is superdiagonal, and $\eta$ is
skew-symmetric. Let $B = e^A$, so that 
$$e^{\mu t} = (B e^{(\delta + \nu)t} B^{-1}) (B e^{\eta t} B^{-1}).$$ 
Since $\eta$ is skew symmetric it follows that
$e^{\eta t}$ is in the orthogonal group $O_n(\R)$.
\end{proof}

It is a surprizing but useful fact that $Ae^BA^{-1} = e^{ABA^{-1}}$.

The theorem concludes that $P^t$ is a bounded subgroup.  That is,
$$\sup_{x\in \R^n} \frac{|P^tx|}{|x|}$$ is bounded for all $t$, where
$|\cdot|$ is the usual Euclidean norm.  This does not 
imply that, using the nilpotent norm $||\cdot||$, $\sup_{x\in
\R^n} \frac{||P^tx||}{||x||}$ is bounded for all $t$, or even that the sup
is finite for a fixed $t$.  For example, consider the Heisenberg group in
$\{x, y, z\}$ basis, with $P \cdot
z = x$.  Then, $\frac{||P(kz)||}{||kz||} = \frac{k}{\sqrt{k}}$.  It is 
precisely this type of map, which fails to preserve the Carnot grading, 
for which $\frac{||Px||}{||x||}$ fails to be bounded. 

In Subsection \ref{subsection:reducetoajf}, we will use Corollary \ref{cor:1pjf}
to show that, with the additional condition that $M$ weakly preserves the
nilpotent grading, this type of mixing is excluded.

\begin{cor}[]\label{cor:1pjf}
Under the conditions of Theorem \ref{theorem:ajf}, 
assume further that matrix $M$ weakly preserves the grading 
$\R^n = V_1 \oplus \cdots \oplus V_k$, with $\dim(V_i) = n_i$.  Then $P^t$ 
is conjugate into the product of the corresponding orthogonal groups 
$O_{n_1}(\R) \times \cdots \times O_{n_k}(\R)$.  
\end{cor}

\begin{proof}
Let $A \in \GL_n(\R)$ be the matrix such that $A^{-1}M A$ is in Jordan 
form.   
Suppose that $V \subset \R^n$ is preserved by $M$. Denote $d = \dim(V)$
and choose a basis $\{e_1, \ldots, e_n\}$ for $\R^n$ such that
$$V = \spn \{e_1, \ldots, e_{d}\}.$$
Express matrix $M$ with respect to this basis.  We define a 
truncated matrix $M' \in \M_{d}(\R)$ as follows:
$$(M')_{i,j} = (M)_{i,j}  \mbox{ if }   1 \leq i,j \leq d.$$
Let $i: V \rightarrow \R^n$ be the inclusion map.  Then, for any $v \in
V$, $i(M' v) = M v$.  Applying Theorem \ref{theorem:ajf} to matrix
$M'$ yields:
$$(M')^t = (\barM')^t (P')^t = (P')^t(\barM')^t,$$
where $\barM', P' \in \M_{d}(\R)$.  In particular, both $\barM'$ and $P'$
preserve $V$.  Since the action of $M$ and $M'$ are identical on $V$, the
action of $P$ and $P'$ must also agree on $V$.  Thus, $P$ must preserve
$V$.  Because $P$ is orthogonal, it must also preserve $V^\perp$.  Thus, 
it acts orthogonally on each component, and $P \in O_d(\R) \oplus 
O_{n-d}(\R)$.  
Applying this argument to each $V_i \oplus \cdots \oplus V_k$ shows
that $P \in O_{n_1}(\R) \times \cdots \times O_{n_k}(\R)$.
\end{proof}

The Jordan form, and thus the absolute Jordan form, of a matrix is unique up
to permutation of the blocks.  When the matrix $M$ represents a linear
transformation of $\R^n$, conjugating represents determining a new choice of
basis for the space $\R^n$, and permuting the blocks corresponds to
permuting the elements of the basis.

\subsection{Reducing a Matrix to Permuted Absolute Jordan Form} 
\label{subsection:reducetoajf}

In this subsection we will complete the proof of Theorem \ref{theorem:pajf}
by first showing that the Lie groups are \qic. This step is mostly handled
by the following theorem, which corresponds to Proposition 4.1 in
\cite{ABC} and is restated here in the setting of \NBC groups.  The proof
requires more subtlety in this situation.

\begin{theorem}[Quasi-isometric Lie Groups]\label{prop:reducetoajf} 
Let $\phi$ and $\theta$ be injective endomorphisms of a fixed nilpotent
Lie group $G$.  Suppose that the maps of the Lie algebra $\g$ induced by
$\phi$ and $\theta$ are represented by matrices $M$ and $N$ which lie on
1-parameter subgroups $M^t, N^t$ in $\GL(n,\R)$. Suppose there exist 
integers $r,s>0$ such that $M^r$ and $N^s$ have the same permuted 
absolute Jordan form. Then the metric spaces $G_\phi$ and $G_{\theta}$ are
quasi-isometric. To be explicit, there exists $A \in \GL(n,\R)$ and $K \ge
1$ such that for each $t \in \R$, the map $v \mapsto A(v)$ is a
$K$-bilipschitz homeomorphism from the metric $d_{\phi,t}$ to the metric
$d_{\theta,\frac{s}{r} \cdot t}$; it follows that the map from $G_\phi =
\R^n \rtimes_{\phi} \R$ to $G_{\theta} = \R^n \rtimes_{\theta} \R$ given
by $$(x,t) \mapsto \left( Ax,\frac{s}{r} \cdot t\right)\ $$ is a
bilipschitz homeomorphism from $G_{\phi}$ to $G_{\theta}$, with
bilipschitz constant $\sup\{K,\frac{s}{r},\frac{r}{s}\}$. \end{theorem}

The proof of this theorem in the setting of \NBC groups will require the 
following two lemmas, which provide information about the effect of $A$ 
and $P^t$ on the nilpotent geometry.

\begin{lemma}[$\phi$ is Bounded.]\label{lemma:A}
Suppose $\g= V_1 \oplus \cdots \oplus V_k$ is the Lie algebra of a Carnot 
group, and that $\phi \colon \g \to \g$ is an injective endomorphism.  
Then $\frac{||\phi(v)||}{||v||}$ is bounded away from 0 and $\infty$.
\end{lemma}

\begin{proof}
Define $\phi_{i,j} \colon V_i \rightarrow V_j$ as follows:
$$\phi_{i,j} (v_i) = \proj_{V_j}(\phi(\inc_i(v_i))),$$
where $\inc_i$ is the inclusion of $V_i$ in $V$.
For each $\phi_{i,j}$, define:
$$m_{i,j} = \min_
{v_i \in V_i \setminus\{0\}}
\left\{\frac{|\phi_{i,j}(v_i)|}{|v_i|}\right\}, 
\mbox{ and } 
M_{i,j} = \max_{v_i \in V_i \setminus\{0\}}
\left\{\frac{|\phi_{i,j}(v_i)|}{|v_i|}\right\}.$$
The minimum and maximum always exist, although they may be zero.
Lemma \ref{theorem:endo} states that $\phi$ weakly preserves the grading 
of 
$\g$ by $\{V_i\}$.  As a result, $\phi_{i,j} = 0$ for $i > j$.  Because 
$\phi$ is also injective and $\g$ is finite dimensional, $m_{i,i} \neq 0$ 
for all $i= 1, \ldots, k$.
\medskip
\noindent
\textbf{Bounded Above:}
Now consider $x = (x_1, \ldots, x_k)$, where each $x_i \in V_i$.  Then,
\begin{equation*}
\begin{split}
\left| \proj_{V_j} \phi(x)\right| &= \left| \Sigma_{1 \leq i \leq j} \phi_{i,j} 
(x_i) \right|\\
&\leq k \cdot \max_{1 \leq i \leq j} \left| \phi_{i,j} (x_i) \right|\\
&\leq k \cdot \max_{1 \leq i \leq j} M_{i,j}|x_i|.\\
\end{split}
\end{equation*}

Corollary \ref{cor:dist2} implies there is $K > 0$ such that:
\begin{equation*}
\begin{split}
||\phi(x)||& \leq K \max_j \left\{\sqrt[j]{(\phi(x))_j}\right\}\\
&\leq K \max_j \left\{k \max_{1 \leq i \leq j} \sqrt[j]{M_{i,j} 
|x_i|}\right\}\\
&\leq K \cdot k \cdot \max_{i,j} \sqrt[j]{\left(M_{i,j}\right)} 
\cdot \max_j \max_{1 \leq i \leq j} \sqrt[j]{|x_i|}.\\
\end{split}
\end{equation*}
Let $M = \max_{i,j} \left\{\sqrt[j~]{M_{i,j}}\right\}$.  Observe also 
that, for fixed $i$,
$$\max_{j \geq i} \sqrt[j]{|x_i|} = 
\sqrt[i]{|x_i|}.$$
Thus, 
\begin{equation}
\label{eqn:bd1}
\begin{split}
||\phi(x)|| &\leq K \cdot k \cdot M \cdot \max_i 
\left\{\sqrt[i]{|x_i|}\right\}\\
&\leq K \cdot k \cdot M \cdot ||x||,
\end{split}
\end{equation}
so $\frac{||\phi(x)||}{||x||}$ is bounded above.

\medskip
\noindent
\textbf{Bounded Below:}
By definition of $m_{i,i}$, 
$$\left|\proj_{V_i}\phi(x)\right| \geq m_{i,i} |x_i|.$$
Define
$$m = \min_i\left\{\sqrt[i]{m_{i,i}}\right\},$$
and note that $m > 0$.  Again, considering the nilpotent metric:
\begin{equation}
\label{eqn:bd2}
\begin{split}
||\phi(x)|| &> \frac{1}{K} \max_i 
\left\{\sqrt[i]{\proj_{V_i}\phi(x)}\right\}\\
&\geq \frac{1}{K} \max_i \left\{ \sqrt[i]{m_{i,i}|x_i|}\right\} \\
&\geq \frac{1}{K}\cdot m \cdot \max_i \left\{ \sqrt[i]{|x_i|}\right\}\\ 
&\geq \frac{1}{K}\cdot m \cdot \frac{||x||}{K}.\\ 
\end{split}
\end{equation}
Thus, $\frac{||\phi(x)||}{||x||}$ is bounded away from both 0 and
$\infty$.
\end{proof}

\begin{lemma}[$P$ is Bounded on $\g$]\label{lemma:rotations} Suppose $\g=
V_1 \oplus \cdots \oplus V_k$ is the Lie algebra of a Carnot group.  
Suppose further that $P \in O(n_1) \oplus \cdots \oplus O(n_k)$, where
$n_i = \dim(V_i)$ and that $P^t$ is bounded on $\R^n$.  Then $$\sup_{x\in
G} \frac{||P^tx||}{||x||}$$ is uniformly bounded; i.e., the bound is
independent of $t$. \end{lemma}

\begin{proof}
Lemma \ref{lemma:A} shows that, for fixed $t$, 
$\sup_{x\in \R^n} \frac{||P^tx||}{||x||}$ is bounded.  In fact, Equations 
\ref{eqn:bd1} and \ref{eqn:bd2} show:
$$\frac{m}{K^2} \leq \frac{||\phi(x)||}{x} \leq MKk,$$
where $K$ and $k$ depend only on the group $G$.  Although $M$ 
and $m$ depend on $P^t$, the fact that $P^t$ is bounded in the Euclidean 
norm implies that $M$ and $m$ are uniformly bounded for all $t$.  This 
implies the lemma.
\end{proof}

\begin{proof}[Proof of Theorem \ref{prop:reducetoajf}]
The proof proceeds similarly to the proof of Proposition 4.1 in 
\cite{ABC} but the nonisotropic nature of the nilpotent geometry makes the 
proof more involved.

\textbf{Case 1:} 
Assume that $N^t = e^{Jt}$ is the unique 1-parameter 
subgroup such that $N= e^J$ is conjugate to the absolute Jordan form of 
$M$.  Then by Theorem \ref{theorem:ajf} and its Corollary \ref{cor:1pjf}
$$M^t = (A^{-1}N^tA)P^t
$$
where $A \in \GL_n(\R)$ and the 1-parameter subgroup $P^t$ is a bounded 
element of the product of the orthogonal groups $O(n_1) \oplus \cdots 
\oplus O(n_k)$. Choose $t \in \R$ and $v \in \R^n$.  We must show that the 
two numbers
$$ ||M^{-t}v|| = ||P^{-1}(A^{-1}N^{t}A)v|| \mbox{  and  } ||N^{-t}Av||
$$
have ratio bounded away from 0 and $\infty$, with bound independent of $t, 
v$. We set $u = N^{-t}Av$, so 
it suffices to show that $||P^{-t}A^{-1}u||$ and $||u||$ have bounded 
ratio.  

As discussed in Subsection \ref{subsection:opjs}, this shall be more difficult
than in the \ABC case.  However, Lemma \ref{lemma:rotations} shows that if
$P \in O(n_1) \oplus \cdots \oplus O(n_k)$ and $P^t$ is bounded, then
$\frac{||P^{t}v||}{||v||}$ remains bounded.  Furthermore, as
shown in Lemma \ref{lemma:A}, the ratio $\frac{||Ax||}{||x||}$ is also
bounded away from both 0 and $\infty$ for the matrix $A$.

\medskip
\noindent
\textbf{Case 2:}
Assume that there exists $a>0$ such that $M^t = N^{at}$ for all  
$t$. Then the metrics $d_{M,t}$ and $d_{N,at}$ are identical.

\medskip
\noindent
\textbf{General case:} 
Applying Case 2 we may assume that $\det M = \det N$. Applying Case 1
twice we may go from $G_M$ to $G_{e^J}$ to $G_N$, where $e^J$ is
conjugate to the absolute Jordan form of $M$ and of $N$.
\end{proof}

We now have the tools to complete the proof of Theorem 
\ref{theorem:pajf}.

\begin{proof}[Proof of Theorem \ref{theorem:pajf}]
Given an endomorphism $\phi$ of a nilpotent group $N$, the \NBC group
$\Gamma_{\phi^k}$ is a finite index subgroup of $\Gamma_\phi$, and thus 
they are \qic.  Now we can restrict our attention to groups defined by 
endomorphisms with the same permuted absolute Jordan form.

Suppose $N_1$ and $N_2$ are lattices in the same Carnot group $G$ and
$\phi_1$ and $\phi_2$ are injective, nonsurjective endomorphisms of $N_1$
and $N_2$ respectively, each acting without unipotent part.  Suppose that
$M_1$ and $M_2$ have the same permuted absolute Jordan form.  By Theorem
\ref{prop:reducetoajf}, the Lie groups associated to $M_1$ and $M_2$ are
quasi-isometric via a height preserving quasi-isometry. It remains to show
there is a height preserving \QI between the associated trees.

If $M'$ is the permuted absolute Jordan form of $M$, then
$$\det(M') = |\det(M)|.$$
Therefore, if $M_1$ and $M_2$ have the same permuted absolute
Jordan form, then their determinants have the same absolute value.  

As described in Subsection \ref{subsection:model}, the tree associated to
$\Gamma_i$ is uniform with 1 `in' and $[N:\phi(N)]$ `out' edges at each
vertex.  The index $[N:\phi(N)] = |det (M)|$, where $M$ is the matrix 
which describes the induced action of $\phi$ on the Lie algebra.  Thus, 
there is a height preserving isometry between the Bass-Serre trees 
associated to $\Gamma_1$ and $\Gamma_2$.

Recall that the group $\Gamma_i$ is \QIc to the fiber product $X_i$ of the
Lie group $G_i$ and the tree $T_i$.  The height preserving
quasi-isometries of the Lie group and the tree induce a \QI of the fiber
product which shows that the \NBC groups $\Gamma_{\phi_1}$ and
$\Gamma_{\phi_2}$ are quasi-isometric. \end{proof}

\subsection{Coarse Topology, inducing a Quasi-isometry of $G_{\phi}$}
\label{subsection:bushy} 
\label{subsection:coarsetop} 
In this subsection, we will establish the following: 

\begin{prop}
\label{prop:coarsetop}
Suppose that $N$ is a discrete nilpotent group which is a lattice in a
nilpotent Lie group $G$.  Let $\phi$ be an injective 
nonsurjective endomophism of $N$.  Recall that $\Gamma_{N, \phi}$ is the
discrete group which is the HNN extension of $N$ by $\phi$.  
Similarly, consider $N'$,
$G'$, $\phi'$, and $\Gamma_{N',\phi'}$. If there exists a
quasi-isometry $f\colon \Gamma_{N, \phi} \to \Gamma_{N', \phi'}$ then:

\begin{enumerate}
\renewcommand{\labelenumi}{\alph{enumi})}
\item There is a quasi-isometry between the nilpotent groups $g \colon N 
\rightarrow N'$.

\item There is a quasi-isometry between the Lie groups $\theta \colon
G_{\phi}\to G'_{\phi'}$ which coarsely respects the transversely oriented
horizontal foliations.  

\item Furthermore, all associated constants for $\theta$ depend only on 
those for $f$.

\end{enumerate}
\end{prop}

This is a modified version of \cite{ABC}, Proposition 7.1.  The proof
there proceeds in four steps.  Our presentation will parallel the one
there, with Steps 1 and 2 (as encapsulated in Theorems 7.3 and 7.7,
respectively) as well as Step 4 applying directly.  Step 3 will require a
new proof, given below in Theorem \ref{theorem:bigon}.  The equivalent of
part (a) in the abelian-by-cyclic case required only showing that
dimension was preserved.  Here, this part is a new Corollary of their
Theorem 7.7.

\begin{proof}
Consider a \QI between two \NBC groups $f \colon \Gamma_{N,\phi}
\rightarrow \Gamma_{N', \phi'}$.  This induces a \QI between the geometric
model spaces $X_\phi$ and $X_{\phi'}$.  
Henceforth, we will think of these spaces and maps interchangably.

\paragraph*{Step 1. Quasi-isometrically embedded hyperplanes are
close to hyperplanes.}

In the language of \cite{ABC}, the HNN extensions $\Gamma_{N, \phi}$ are
finite, geometrically homogenous graphs of groups.  Furthermore, the edge and
vertex groups are fundamental groups of Poincare duality spaces.  This is 
because (the universal cover) $G$ is connected and simply-connected, and 
therefore homeomorphic to $\R^n$.  Also, $\pi_1(G/\Gamma) = \Gamma$.  
Since the manifolds $G/\Gamma$ satisfy
Poincare duality, each $\Gamma$ is a Poincare duality group.  Thus, as
discussed in \cite{ABC}, the following theorem applies to the metric 
fibration of the geometric model space over the associated tree: $X_{N, 
\phi} \rightarrow T_{N, \phi}$.

\begin{theorem}[\cite{ABC}, Theorem 7.3]
Let $\pi \colon X \rightarrow T$ be a metric fibration whose fibers are
contractible $n$-manifolds for some $n$. Let $P$ be a contractible
$(n+1)$-manifold which is a uniformly contractible, bounded geometry,
metric simplicial complex. Then for any uniformly proper embedding $\phi
\colon P \rightarrow X$, there exists a unique hyperplane $Q\subset X$ 
such that 
$\phi(P)$ and $Q$ have finite Hausdorff distance in $X$. The bound on
Hausdorff distance depends only on the metric fibration data for $\pi$,  
the uniform contractibility data and bounded geometry data for $P$, and
the uniform properness data for~$\phi$.
\label{theorem:metricfibration}
\end{theorem}

So, each hyperplane of $X_{N,\phi}$ is mapped by $f$ to a (universally)
bounded neighborhood of some hyperplane $P_{\ell'} \subset X_{N',\phi}$.

\paragraph*{Step 2. A quasi-isometry takes hyperplanes and horizontal  
leaves in $X_\phi$ to hyperplanes and horizontal leaves in $X_{\phi'}$:} 

This is the step in the proof which depends upon having endomorphisms of
the discrete groups $N$, $N'$ which are \emph{not} surjective.  This
implies that the tree $T= T_{N,\phi}$ (respectively, $T' = T_{N',\phi'}$)
has uniform valence $[N:  \phi(N)] > 1$(respectively, $[N':  \phi'(N')] >
1$).  For $\beta > 0$ we say a tree $T$ is \emph{$\beta$-bushy} if each 
point of $T$ is within distance
$\beta$ of some vertex $v$ such that $T-\{v\}$ has at least 3 unbounded
components.  So, the trees $T$ and $T'$ are $\tfrac{1}{2}$-bushy, in the
path metric described in Subsection \ref{subsection:model}.

Thus, the following theorem applies:
\begin{theorem}[\cite{ABC}, Theorem 7.7]\label{theorem:horizontal}
Let $\pi \colon X \to T$, $\pi' \colon X' \to T'$ be  metric fibrations
over $\beta$-bushy trees $T,T'$, such that the fibers of $\pi$ and $\pi'$
are contractible $n$-manifolds for some $n$. Let $f\colon X \to X'$ be a
quasi-isometry. Then there exists a constant $A$, depending only on the
metric fibration data of $\pi,\pi'$, the quasi-isometry data for $f$, and
the constant $\beta$, such that:
\begin{enumerate}\renewcommand{\labelenumi}{\alph{enumi})}
\item For each hyperplane $P \subset X$ there exists a unique
hyperplane $Q \subset X'$ such that $d_\Haus(f(P),Q) \le A$.
\item For each horizontal leaf $L \subset
X$ there is a horizontal leaf $L' \subset X'$ such that 
$d_\Haus(f(L),L')\le A$.
\end{enumerate}
\end{theorem}

As an immediate consequence, we have:

\begin{cor}[Quasi-isometric Base Groups]
\label{cor:basegp}
If $\Gamma_{N,\phi}$  and  $\Gamma_{N',\phi'}$ are \NBC groups 
which are quasi-isometric, then $N$ and $N'$ are quasi-isometric, and 
there is a quasi-isometry 
$f: G_{N,\phi} \rightarrow G_{N',\phi'}$ which coarsely respects 
the horizontal 
foliations and their transverse orientiations.
\end{cor}

\begin{proof} 
Let $X = X_{N,\phi}$ be the geometric model space for $\Gamma_{N,\phi}$, 
and $X'$ be the model space for $\Gamma_{N',\phi'}$.
If $f \colon X \rightarrow X'$ is a $(K,C)$-\qi, then Theorem 
\ref{theorem:horizontal} implies that
there is $A > 0$ such that for every horizontal leaf $L \subset X$, 
there is a horizontal leaf $L'
\subset X'$ such that $d_\Haus(f(L),L')\le A$.  Let $p \colon X'
\rightarrow L'$ be nearest point projection.  Then, $p \circ f \colon L
\rightarrow L'$ is a $(K, C + 2A)$-\qi.  Furthermore, $N_{2A}(p \circ
f(L)) \supset L'$.  Thus $p \circ f$ shows that $L$ and $L'$ are
quasi-isometric.  Since $L$ is quasi-isometric to the nilpotent discrete
group $N$, and $L'$ is quasi-isometric to the nilpotent discrete group
$N'$, we conclude that $N$ and $N'$ are quasi-isometric. \end{proof}

This result shows that, in classifying \NBC groups, we may restrict our
attention to groups which are extensions of nilpotent groups which are
quasi-isometric.  As mentioned in the introduction, the quasi-isometric
classification of nilpotent groups is a major unsolved problem.  For the 
classification given here, we will make the slightly stronger

{\bf Standing Assumption:} Discrete nilpotent groups $N$ and $N'$ are 
lattices in the same Carnot group $G$.

\paragraph*{Step 3. A quasi-isometry takes coherent hyperplanes in     
$X_\phi$ to coherent hyperplanes in $X_{\phi'}$.}

\begin{theorem}[Coherence is a \QI invariant]\label{theorem:bigon}
Given a \QI between hyperplanes $f:H_1 \rightarrow H_2$, the hyperplane 
$H_1$ is coherent if and only if $H_2$ is coherent.
\end{theorem}

As in the \ABC case, the proof of this step is based on a comparison of 
the growth functions
for coherent and incoherent hyperplanes.  For coherent hyperplanes, the
growth will be linear (it is quadratic in the \ABC case), while for
incoherent hyperplanes it will be exponential. 

\begin{proof}
We need the following definitions as in \cite{ABC}:
For any hyperplane $H \subset X$, there is a quotient map $H \to \R$
whose point pre-images give the horizontal foliation of $H$, and such
that the Hausdorff distance between two horizontal leaves equals the
distance between the  corresponding points in $\R$. A path $\gamma$ in
$H$ is said to be  $(K,C)$-\emph{quasivertical} if its projection to $\R$
is a $(K,C)$-quasigeodesic. Define a $(K,C)$-\emph{quasivertical bigon}
in $H$ to  be a pair of $(K,C)$-quasivertical paths $\gamma,\gamma'$
which begin and  end at the same point. 

If $K,C$ are fixed, we define a filling area function $A(L)$ for
$(K,C)$-quasivertical bigons in $H$. Given a $(K,C)$-quasivertical
bigon $\gamma,\gamma'$, its \emph{filling area} is the infimal
area of a Lipschitz map $D^2 \to H$ whose boundary is a
reparameterization of the closed curve $\gamma^{-1}*
\gamma'$; such a map $D^2 \to H$ is called a \emph{filling
disc} for $\gamma^{-1}* \gamma'$. For each $L \ge 0$ define
$\A(L)$ to be the supremal filling area over all $(K,C)$-quasivertical
bigons $\gamma,\gamma'$ in $H$ such that $\Length(\gamma) +
\Length(\gamma') \le L$.

Suppose that two hyperplanes $H_1$ and $H_2$ have filling functions
$\A_1(L)$ and $\A_2(L)$, respectively.  As is shown in \cite{ABC}, if
there is a \QI between $H_1$ and $H_2$, then the filling functions must be
comparable in the following sense: \begin{equation}\label{equation:bigon}
\A_1(L) \leq \alpha \cdot \A_2(\beta L + \delta) + \zeta, \end{equation}
for some $\alpha, \beta, \delta, \zeta > 0$ which do not depend on $L$, 
and the same equation must hold with $\A_1$ and $\A_2$ reversed.

However, we will see that in the geometric model space associated to a
\NBC group, the filling function of a coherent hyperplane is linear, while
the filling function of an incoherent hyperplane is exponential.

Consider a $(K,C)$-quasivertical bigon in a coherent hyperplane.  As shown
in Subsection \ref{subsection:flow}, a quasivertical line must be shadowed by a
vertical flow line.  (This is a result of the assumption that the
endomorphism is unipotent-free; without this assumption, we could get
quadratic growth, as in the \ABC case.)  Thus, the distance between the
two edges of the bigon is universally bounded.  Therefore, the filling
area is linear in $L$.

In the case of an incoherent hyperplane, the argument used in \cite{ABC}
depends only on having one direction in which the growth of vectors is
exponential with base greater than 1.  This is true in the \NBC situation,
although the growth is calculated as $\lambda^{\frac{1}{w_i}}$ and not
simply $\lambda$.

Since linear and exponential functions are not comparable according to
equation \ref{equation:bigon}, there must not be a \QI between coherent 
and
incoherent hyperplanes. 
\end{proof}

\paragraph*{Step 4. A horizontal-respecting quasi-isometry preserves
transverse orientation.}

The proof in \cite{ABC} applies directly.  The idea of the proof is this:
a quasi-isometry must preserve the sign of the log of the determinant.  

That is, given a \QI $\phi \colon \Gamma_1 \rightarrow \Gamma_2$, the
absolute values of the determinants of the matrices associated to 
$\Gamma_1$ and $\Gamma_2$ must
either both be greater than one or both less than one.  Since reversing
orientation comes from taking the inverse of the matrix, this shows that
the \QI must preserve orientation.

This concludes the proof of Proposition \ref{prop:coarsetop}.
\end{proof}

\subsection{Time Change Rigidity, Part 1}\label{subsection:time1}

Recall that a map $f$ between hyperplanes $P$ and $P'$ is 
\emph{horizontal-respecting} if there exists a function $h \colon \R
\rightarrow \R$ and $A \ge 0$ such that
$d_\Haus(f(P^{\vphantom\prime}_{t}),P'_{h(t)}) \le A$ for all $t\in \R$, 
where $P_t = \pi^{-1}(\ell \cap \{t\})$.
The function $h \colon \R \rightarrow \R$ is called an \emph{induced time
change} for $f$ with \emph{Hausdorff constant} $A$.  
If $h$ and $h'$ are
two induced time changes for $f$ with Hausdorff constants $A$ and $A'$,
then $\sup_t |h(t) - h'(t)| \leq A + A'$.
The converse is also true: If  $h$ is an induced time change for $f$ 
with Hausdorff constant $A$, and $h': \R \rightarrow \R$ satisfies
$\sup_t |h(t) - h'(t)| \leq A'$, then $h'$ is also an induced time change 
function for $f$ with Hausdorff constant $A + A'$.

Let $f$ be the \QI between hyperplanes which is guaranteed by Corollary
\ref{prop:coarsetop}. Later (Subsection \ref{subsection:time2}), we will show
that there is a linear function which is an induced time change function
for $f$.  That proof will depend upon the result of Subsection
\ref{subsection:flow}: that vertical flow lines are coarsely preserved.  To
show that, however, we need the fact that the induced time change is at
least coarsely linear and coarsely increasing; that is, $h \colon \R
\rightarrow \R$ is a \QI and $h(t) \rightarrow \infty$ as $t\rightarrow
\infty$.  The proof of the following Lemma, given by Farb and Mosher
(\cite{ABC}, Lemma 5.1) for the \ABC case, applies in this context as
well.

\begin{Lemma}\label{lemma:timeqi} For each $K,C,A$ there exists $C'$ 
such
that if $f \colon G_\phi \rightarrow G_{\phi'}$ is a horizontal respecting
$(K,C)$ quasi-isometry, and $h \colon \R \rightarrow \R$ is an induced
time change for $f$ with Hausdorff constant $A$, then $h$ is a $(K,C')$
quasi-isometry of $\R$. \end{Lemma}

\subsection{Vertical Flow Lines are Coarsely Preserved}\label{subsection:flow}
In this subsection, we will show that the \QI $f$ between $G_{N, \phi}$ and
$G_{N',\phi'}$ established by Corollary \ref{prop:coarsetop} coarsely
preserves vertical flow lines.  As explained in \cite{ABC}, if all the
eigenvalues of $M$ and $N$, the matrices associated to $\phi$ and $\phi'$,
respectively, are greater than 1, then $G_M$ and $G_N$ are negatively
curved, and this is simply the fact that a quasigeodesic in a negatively
curved space is Hasudorff close to a geodesic.  Since we permit $M$
to have any eigenvalue with absolute value different from 1, we
need the stronger:

\begin{theorem}\label{theorem:flow}
Consider a \QI $f \colon G_\phi \rightarrow G_{\phi'}$.  There exists 
$\alpha \geq 0$ such that 
for each vertical flow line $l_x$ in $G_\phi$, there exists
a vertical flow line $m_y$ in $G_{\phi'}$ such that $f(l_x)$ is
contained in the $\alpha$-neighborhood of $m_y$.
\end{theorem}

The proof given in \cite{ABC} (Claim 5.7, p.165-6) goes through exactly as
given.  The idea of the proof is given below, following some definitions.  
In that context, the theorem is stated for leaves of the center foliation,
but with our additional assumption that $\Gamma$ has no center (implied by
the standing assumption that $\phi$ is unipotent-free), the center leaves
are simply vertical flow lines.

\begin{Definition}
Given a flow $\Phi$ on a metric space $X$, we write $x\cdot t$ as an
abbreviation for $\Phi_t(x)$.  Given $\epsilon, T \geq 0$, an $(\epsilon,
T)$-pseudo-orbit of $\Phi$ consists of a sequence of flow segments $(x_i
\cdot [0, t_i])$, where the index $i$ runs over an interval in $\Z$, such
that $d_X(x_i \cdot t_i, x_{i+1}) \leq \epsilon$ and $t_i \geq T$ for all
$i$. \end{Definition}

A flow $\Phi$ on a manifold $M$ is called \emph{hyperbolic} if, at 
each $x \in M$, the tangent space has a splitting $T_xM=E^u \oplus 
E^x$ such that the derivative $d\Phi: TM \rightarrow TM$ expands $E^u$, 
contracts $E^s$, and preserves both.  The subbundles $E^s$ and $E^u$ are 
tangent to the global stable and unstable manifolds, defined by:
$$W^s(x) = \{y \in M \suchthat d(f^n(x), f^n(y)) \rightarrow 0 \mbox{ as } 
n \rightarrow \infty\}, \mbox{ and }$$
$$W^u(x) = \{y \in M \suchthat d(f^{-n}(x), f^{-n}(y)) \rightarrow 0 
\mbox{ as } 
n \rightarrow \infty\}.$$

Given these structures, 
it can be shown (e.g., \cite{KH} Proposition 6.4.13) that:
\begin{itemize}
\item there is an 
$\epsilon > 0$ such that for any $x, y \in X$, the intersubsection 
$W^s(x) \cap N_\epsilon (x) \cap W^u(y) \cap N_\epsilon(y)$ consists of at 
most one point, and  
\item there is a $\delta >0$ such that $W^s(x) 
\cap N_\epsilon (x) \cap W^u(y) \cap N_\epsilon(y) \neq \emptyset$ 
whenever $d(x,y) < \delta$.
\end{itemize}

Such a structure is called a \emph{local product structure}.  
An equivalent statement is: there is a $\delta >0$ and, for all $p \in X$, 
a map $f_p \colon
N_\delta(p) \rightarrow 
\R^k \times \R^{n-k}$ such that $f(W^u(x)) \subset \R^k
\times \{v\}$ for some $v \in \R^{n-k}$ and $f(W^s(x)) \subset \{w\}
\times \R^{n-k}$ for some $w \in \R^k$.

If there is a single map $f: M \rightarrow \R^k \times \R^{n-k}$ which
satisfies the above conditions, for $ \epsilon, \delta = \infty$ then we
say $f$ has a \emph{global product structure}.

The existance of a local product structure for hyperbolic flows is a key
ingredient in the proof of the so-called Shadowing Lemma.  This was first
proved by Bowen (\cite{Bowen}, Theorem 2.2, Approximation Theorem).  
Hirsch, Pugh, and Shub found a new proof for the following more 
streamlined statement:

\begin{Lemma}[Shadowing Lemma, \cite{HPS} Lemma 7A.2, p. 133] If $(f, L)$
has local product structure and $\nu > 0$ is given, then there exists
$\delta > 0$ such that any $\delta$-pseudo orbit for $f$ in $\Lambda$ can
be $\nu$-shadowed by a pesuo-orbit for $f$ in which $\Lambda$ which
respects $L$.  \end{Lemma}

Notice the quantifiers: in this statement (as in Bowen's presentation),
$\delta$ depends on $\nu$, and as $\nu$ approaches zero, so does $\delta$.  

In Bowen's proof, dependence of $\delta$ on $\nu$ ($\epsilon$ in his
statement) arises in two different ways.  First, he needs local product
structures on neighborhoods with radius at least $\delta$.
This dependence is not surprizing; it says that we must have a product
structure on neighborhoods which are large enough that the pseudo-orbit 
cannot jump out of the neighborhood. Secondly, he
requires that $\delta$ be small enough to satisfy a bound on the sum of a
particular convergent geometric sequence.  This sequence also depends upon 
$\min\{\lambda >1, \lambda^{-1} >1 \suchthat \lambda \mbox { is an
eigenvalue of } $df$\}$.  This condition can be satisfied for a fixed 
$\delta$ by raising $M$ to large enough powers.  

In applying this lemma here (as in \cite{ABC}), we want to choose an
arbitrarily large $\delta$.  This can be accomplished if the requirement
of a local product structure is replaced by the condition of having a
global product structure. Farb and Mosher saw that the existance of a
global product structure is sufficient to ensure that, for arbitrary
choice of $\delta$, there is a $\nu$ so that the conclusion of the
Shadowing Lemma holds.  They stated and used the following ``Global
Shadowing Lemma''. (It is stated here in the special case that center
leaves are vertical flow lines, which is true for \NBC groups.)

\begin{Lemma}[Shadowing Lemma, \cite{ABC} Lemma 5.3, p. 163] Consider a
1-parameter subgroup $M^t$ of $GL_n(\R)$, and let $\Phi$ be the vertical
flow on $G_M$.  For every $\epsilon, T \geq 0$, there exists $\delta \geq
0$ such that every $(\epsilon, T)$-pseudo-orbit of $\Phi$ is
$\delta$-shadowed by a vertical flow line $m_y$.  That is, if $(x_i \cdot
[0, t_i])$ is an $(\epsilon, T)$-pseudo-orbit, then there is a vertical
line $m_y$ such that $d(x_i \cdot t, y \cdot t) \leq \delta$ for all $i$
and all $t \in [0,t_i]$. \end{Lemma}

In the cases of \ABC and \NBC groups, a global product structure exists
for any hyperbolic flow on $G$.  On the Lie algebra $\g$ of $G$, the
global stable manifold $W^s$ is just the span of the root spaces with
eigenvalues less than 1;  the global unstable manifold is the span of the
eigenspaces greater than 1.  Here, the global product structure is clear; 
the leaves are all linear subspaces, and $\dim(W^s) + \dim(W^u) = n$, so 
any stable manifold intersects any unstable manifold at a unique point.
 
The exp map is a diffeomorphism from $\g$ to $G$ which conjugates $df$ to
$f$.  Thus, it carries stable manifolds of $\g$ to stable manifolds of $G$
and preserves the global product structure.  Thus, $G$ has a global
product structure as well.

Under our condition of unipotent-free endomorphisms, this gives a global
decomposition of the space. Thus, Farb and Mosher's Global Shadowing Lemma
applies in the situation of \NBC groups.

\begin{proof}[Idea of Proof of Theorem \ref{theorem:flow}] The proof is in
two steps.  First, show that $f(l_x)$ is close to a pseudo-orbit.  To do
this, use the fact that there is an coarsely-linear induced time change
which is coarsely increasing (Lemma \ref{lemma:timeqi}) to choose a sparce
but regularly spaced sequence $\{x_i\}$ in the flow $l_x$.  The images
$\{y_i=f(x_i)\}$ of these points are used to define a pseudo-orbit in
$G_{\phi'}$.  Then, the Shadowing Lemma implies that the pseudo-orbit is
within a bounded neighborhood of a vertical flow line.
\end{proof}

\subsection{Time Change Rigidity, Part 2}\label{subsection:time2}

In Subsection \ref{subsection:time1}, we showed that the induced time change was 
coarsely linear and coarsely increasing.  Now, with the additional 
information that vertical flow lines are coarsely preserved, we can show 
more:

\begin{theorem}[Linear Induced Time Change]\label{theorem:time} Suppose
$f$ is a horizontal-respecting \QI between Lie groups $G_\phi$ and
$G_{\phi'}$, both of which satisfy our standing hypotheses. (That is,
$\phi$ and $\phi'$ are unipotent-free nonsurjective endomorphisms of
discrete nilpotent groups $N$ which are lattices in the same Lie group
$G$.)  Then, there exists $m > 0$ such that $h(t) = mt$ is an induced time
change for $f$. \end{theorem}

The proof relies upon comparing the divergence rates of vertical flow
lines in each Lie group. The corresponding result in the \ABC case
(\cite{ABC}, Prop. 5.8)  is similar but can be proved much more simply.  
In that case, the smallest divergence rate which is greater than the
constant must be a pure exponential.  That is, it must have the form
$\lambda^t$.  In that case, if $\alpha$ is the eigenvalue of $G_\phi$ with
minimal absolute value greater than 1, and $\beta$ is the corresponding
eigenvalue for $G_{\phi'}$, then the slope of the induced time change is
$m = \frac{\log \alpha}{\log \beta}$.

Much more work will be required in the \NBC case.  It seems true, but
perhaps difficult to show, that the smallest growth rate is pure
exponential.  (Perhaps the smallest growth rate is
$(t^k\lambda^t)^{\frac{1}{w}}$ for some $k \neq 0$, while the vectors
which grow as $\lambda^t$ are in a grade of the nilpotent group with
smaller weight $w$.)  Instead, we compute the induced time change which
would be implied if a line with pure exponential growth were taken instead
to a polynomial$\cdot$exponential diverging line.  Considering the coarse
inverse of $f$, we find a similar rate for the inverse.  Putting these
together shows that $h$ must be linear.

The following Lemma establishes induced time change parameters, given the 
divergence rates of flow lines which are preserved by the \qi.

\begin{lemma}[]\label{lemma:time} 
Consider a horizontal respecting,
vertical flow line preserving $(K,C)$-\QI $f: G_\phi \rightarrow
G_{\phi'}$ for $K \geq 1, C \geq 0$.  Let $g$ be the map between the 
vertical line spaces, (with
constant $R$)  and suppose that $g(\ell_1) = m_1$ and $g(\ell_2) = m_2$,
where $\ell_i$ is a vertical flow line in $G_{\phi}$ and $m_i$ is a
vertical flow line in $G_{\phi'}$.  Suppose further that $d(\ell_1(t),
\ell_2(t))\simeq t^k \cdot \mu^t$ and $d(m_1(t), m_2(t))\simeq \lambda^t$ 
with $k >0$ and $\lambda, \mu > 1$.  
Let $K_\ell \geq 1$ and $C_\ell \geq 0$ be the constants implicit in the 
equation
$d(\ell_1(t), \ell_2(t))\simeq t^k \cdot \mu^t$, and define $K_m, 
C_m$ similarly.  Then, there are positive 
constants $m, c$ such that the
function $h(t) = mt + c\log (t)$ is an induced time change for $f$.
\end{lemma}

\begin{proof} 
We will find two 
different bounds for the distance $d(f(\ell_1(t)), 
f(\ell_2(t)))$, which, when compared, give restrictions on the induced 
time change function.  Suppose that $h(t)$ is an induced time change 
function for $f$ with Hausdorff constant $A$.  (One must exist, since $f$ 
is horizontal foliation 
preserving.)  Then,
$$d(f(\ell_i(t)), m_i(h(t))) < R + A.$$
Consequently,
$$ |d(m_1(h(t)), m_2(h(t))) - d(f(\ell_1(t)), f(\ell_2(t)))| \leq  2R + 
2A,$$
and applying the known divergence rate of of $m_i$ yields:
\begin{equation}\label{1}
\frac{1}{K_m} \mu^{h(t)} - C_m - 2R - 2A \leq
d(f(\ell_1(t)), f(\ell_2(t))) \leq 
K_m \mu^{h(t)} + C_m +2R + 2A.
\end{equation}

On the other hand, by definition of \qi, we have 
$$\frac{1}{K} \cdot d(\ell_1(t), \ell_2(t)) - C \leq d(f(\ell_1(t)),
f(\ell_2(t))) \leq K \cdot d(\ell_1(t), \ell_2(t)) + C,$$
and substituting the divergence rate of $\ell_i$ yields:
\begin{equation}\label{2}
\frac{1}{KK_\ell} \cdot t^k\lambda^t - \frac{C_\ell}{K} - C \leq 
d(f(\ell_1(t)),f(\ell_2(t))) \leq 
KK_\ell \cdot t^k\lambda^t + KC_\ell + C.
\end{equation}

Thus, the upper bound of Equation \ref{1} must be greater than the lower
bound of Equation \ref{2}, and vice versa.  That is, $$K_m \mu^{h(t)} +
C_m +2R + 2A \geq \frac{1}{KK_\ell} \cdot t^k\lambda^t -\frac{C_\ell}{K} - 
C
$$ and $$KK_\ell \cdot t^k\lambda^t + K C_\ell + C \geq \frac{1}{K_m}
\mu^{h(t)} - C_m - 2R - 2A.$$
Using the new constants $K' = K K_\ell K_m$ and $C' = C + C_\ell + C_m + 
2R +2A$, we 
combine these two equations to get:
$$\frac{1}{K'} \cdot t^k\lambda^t - C' \leq \mu^{h(t)} \leq 
K' \cdot t^k\lambda^t + C'.$$

Notice that, for fixed $C' > 0$, and $x \geq C'$, we have 
$$\log(x + C') \leq \log x + \log 2.$$  
When $x < C'$, we have $$\log(x + C') < \log(2C').$$  
Let $M = \max\{\log 2, \log(2C')\}$, so for all $x$, 
$$\log(x + C') \leq 
\log x + M.$$  
Similarly, for $x > C'$, we have 
$$\log (x - C') \leq \log x - M,$$ 
where $M = \max \{\log 2, log(2|C'|\}$.
Taking logs yields: 
$$ -\log K' + k\log t +t \log \lambda - M \leq h(t) \cdot
\log\mu \leq \log K' +k \log t + t \log \lambda + M.$$ 
Thus, 
$$ |h(t) - \frac{k \log t + t \log \lambda}{\log \mu}| \leq \frac{\log 
K' + M}{\log \mu}.$$

This shows that 
$$h'(t) =  \frac{k \log t + t \log \lambda}{\log \mu}$$ 
is an induced time change function for $f$ with Hausdorff constant
$\frac{\log K' + M}{log \mu}$. \end{proof}

Now, we will prove Theorem \ref{theorem:time}.
\begin{proof} 

From Theorem \ref{theorem:divrates} (c), we know that $G_\phi$ has at
least one vertical flow line $\ell$ with a purely exponential divergence
rate, say, $\lambda^t$.  By Theorem \ref{theorem:flow}, $f(\ell)$ is close
to some vertical flow line, say $m$ in $G_\phi$.  By Theorem
\ref{theorem:divrates} (b), the divergence rate of $m$ is in $\C$, that
is, it is of the form $t^{\frac{k}{w}} \cdot (\nu^{\frac{1}{w}})^t$. By
Lemma \ref{lemma:time}, there exist $c, m$ such that $f$ has induced time 
change function $$h(t) = c \log t + m t.$$ Similar reasoning applies to 
the coarse inverse of $f$
(denoted $f'$) to show that there are $c', m'$ such that $f'$ has an 
induced time change $$h'(t) =
c'\log t + m't.$$
Since $f \circ f'$ is a bounded distance from the identity, $h \circ h'$
should also be a bounded distance from the identity.  Yet, both
$\frac{h(t)}{t}, \frac{h'(t)}{t} \rightarrow \infty$ as $t \rightarrow
\infty$, if $k > 0$ or $k' > 0$, respectively.  Thus, $\frac{h(h'(t))}{t}
\rightarrow \infty$ unless both $c, c' = 0$. This occurs only 
when $k, k' = 0$.
\end{proof}

\subsection{Growth Spaces}\label{subsection:growth}

\begin{definition}
Given a Lie group $G_\phi$ and $\lambda > 0$, define the 
\emph{$\lambda$ growth subspace} as
$$\g_\lambda = \{v \in \g \suchthat ||M^tv|| \preceq \lambda^t t^k\mbox{ 
for some 
} k \in \N \}.$$
\end{definition}

This subsection will be devoted to proving:

\begin{prop}[Growth Spaces are Preserved]\label{prop:growthspaces}
If $f: G_\phi \to G_\theta$ is a \QI which preserves the growth space, 
then for all $\lambda \in \R^\times$ 
\end{prop}

In order to prove this, we will need the following

\begin{lemma}[Growth Spaces are Subalgebras]\label{lemma:subalg}
For all $G_\phi$ and $\lambda \in \R^{\times}$, the $\lambda$-growth 
space $\g_\lambda \subset \g$ is a subalgebra.
\end{lemma}

\begin{proof}
By Proposition \ref{theorem:linalg}, any $M$-invariant subspace is 
$\spn\{e_1, \ldots, e_k\}$ for some $k$.  
Suppose $e_i, e_j \in \g_\lambda$.  It suffices to prove $[e_i, e_j] \in 
\g_\lambda$.

Lemma \ref{lemma:divrates1} describes precisely the growth of each 
component of $M^te_i$.    In the notation given there:

\begin{equation*}
\begin{split}
M^t [e_i, e_j] &= [M^te_i, M^te_j] \\
&= [\Sigma_{l=0}^L \lambda_i^{t-l}t^l 
e_{i+l}, \Sigma_{l'=0}^{L'} \lambda_j^{t-l'}t^{l'} e_{j+l'}]\\
&= \lambda_i^t \lambda_j^t \Sigma_{l=0}^L \Sigma_{l'=0}^{L'} 
\lambda_i^{-l}
\lambda_j^{-l'}t^{l+l'}[e_{i+l}, e_{j+l'}]\\
&\simeq (\lambda_i \lambda_j)^t \Sigma_{l=0}^L \Sigma_{l'=0}^{L'} 
t^{l+l'}[e_{i+l}, e_{j+l'}].\\
\end{split}
\end{equation*}

The constant implicit in the transition to $\simeq$ is $K = \max 
\{\lambda_i^{\pm L} \cdot \lambda_j^{\pm L'}\}$.

Proposition \ref{theorem:linalg} (b) implies that, within a Jordan block, 
weights are nondecresing.  That is,
$w_{i+l} \geq w_i$, 
and 
$w_{j+l'} \geq w_j$, so 
$w_{i+l} + w_{j+l'} \geq w_i + w_j$.  Since $G$ is Carnot, $[e_{i + k}, 
e_{j + k}] \in V_{w_{i+l} + w_{j+l'}}$, so 
$$\Sigma_{l=0}^L \Sigma_{l'=0}^{L'} t^{l+l'}[e_{i+l}, e_{j+l'}] \in V_{w_i 
+ w_j} \oplus \cdots \oplus V_k.$$
Furthermore, when this summation is written as a linear combination of the 
basis elements, the coefficients will all be polynomials in $t$ with 
degree at most $l + l'$.  That is, there are real numbers $\alpha_k$ and 
real polynomials $p_l(t)$ such that
$$\Sigma_{l=0}^L \Sigma_{l'=0}^{L'} t^{l+l'}[e_{i+l}, e_{j+l'}] = 
\Sigma_{k=0}^n \alpha_k p_k(t)e_k.$$
Now, apply Corollary \ref{cor:dist1}
to find the nilpotent length of these vectors:
$$||M^t[e_i, e_j]|| \simeq \max_l \{((\lambda_i \lambda_j)^t \alpha_l 
p_l(t) )^{\frac{1}{w_l}}\}.$$

Let $L$ be the value of $l$ which accomplishes the maximum above.  Then

\begin{equation}\label{equation:compare}
\begin{split}
||M^t[e_i, e_j]|| 
&\simeq ((\lambda_i \lambda_j)^t p_L(t) )^{\frac{1}{w_L}}\\
&\simeq ((\lambda_i \lambda_j)^{\frac{1}{w_L}})^t(p_L(t))^{\frac{1}{w_L}}.
\end{split}
\end{equation}

Now turn we direct our attention to the growth rates of $e_i$ and $e_j$.  
Assume, without loss of generality, that 
$$(\lambda_i)^{\frac{1}{w_i}}\geq (\lambda_j)^{\frac{1}{w_j}}.  $$
Then,
\begin{equation*}
\begin{split}
\left(\lambda_i\right)^{\frac{w_j}{w_i}}&\geq \lambda_j\\
\left(\lambda_i\right)^{\frac{w_i + w_j}{w_i}}&\geq \lambda_j\lambda_i\\
\left(\left(\lambda_i\right)^{\frac{1 }{w_i}}\right)^{\frac{ w_i + w_j 
}{w_k}}&\geq 
\left(\lambda_j\lambda_i\right)^{\frac{1}{w_k}}.
\end{split}
\end{equation*}
Since $\frac{ w_i + w_j }{w_k} \leq 1$, we get that
$$(\lambda_i)^{\frac{1 }{w_i}} \geq \left(\left(\lambda_i\right)^{\frac{1 
}{w_i}}\right)^{\frac{w_i + w_j }{w_k}} \geq 
\left(\lambda_j\lambda_i\right)^{\frac{1}{w_k}}.$$

By assumption, $\lambda \geq \lambda_i^{\frac{1}{w_i}}$, so $[e_i, e_j] 
\in \g_\lambda$.
\end{proof}

\begin{proof}[Proof of  Proposition \ref{prop:growthspaces}]

By Lemma \ref{lemma:subalg}, the growth subspace $\g_\lambda$ is a 
subalgebra.  Thus, $exp(\g_\lambda)$ is a connected, simply-connected 
nilpotent Lie group.
A powerful theorem of Pansu (\cite{P}, Theorem 3 in Subsection 1) states 
that if two connected, simply-connected nilpotent Lie groups are
quasi-isometric, then the associated graded Lie algebras are isomorphic.  
In particular, the dimension of each grade is the same. \end{proof}

\section{Some Illustrative Examples}
\label{section:examples}
In each subsection of this section, we present a pair of nilpotent 
groups 
which illustrate the potential complexity in classifying such groups via 
their divergence rates and growth spaces.  The first pair of groups have the 
same permuted absolute Jordan form (and therefore are \qic) although the \QI 
between them is far from being a homomorphism.  The second pair of groups 
agree on the \QI invariants found in this paper but differ in permuted 
absolute Jordan form.  It remains to be determined whether they are \qic.

Let $N = H \times H \times H$, where $H$ is the Heisenberg group.  Thus, 
$N$ is a two-step nilpotent group with presentation:
$$N = \langle a_1, a_2, a_3, a_4, a_5, a_6, a_7, a_8, a_9 \suchthat [a_1, a_2] 
= a_3, [a_4, a_5] = a_6, [a_7, a_8] = a_9 \rangle.$$

All of the endomorphisms here satisfy the property: If $G < N$ is the
infinite cyclic subgroup generated by any one of the nine generators, then
$\phi$ preserves $G$.  That is, for each $a_i$ there is an integer $k_i$ such 
that
$\phi(a_i) = a_i^{k_i}$.  Therefore, we can describe the function $\phi$
by specifying the 9-tuple of exponents.  In fact, all of our $k_i$ will be 
powers of 2, so we prefer to
keep track of $n_i = \frac{1}{2}\log_2 k_i$ for $i=3, 6, 9$ and 
$n_i = \log_2 k_i$ for $i\neq 3, 6, 9$.  
The relations of $N$ imply: $n_3 = \frac{n_1 + n_2}{2}$, $n_6
= \frac{n_4 + n_5}{2}$, and $n_9 = \frac{n_7 + n_8}{2}$.

The matrix $M$ representing such an endomorphism is the diagonal matrix with 
$M_{i,i} = k_i$
Thus, the permuted absolute Jordan 
form of the matrix is 
found by permuting the diagonal entries so that $k_3, k_6$, and $k_9$ 
are the entries in the upper left, in increasing order, followed by $\{k_1, 
k_2, k_4, k_5, k_7, k_8\}$, also in increasing order.

\subsection{Quasi-isometric Groups Can Be Quite Different}
\label{subsection:example1}

Compare the endomorphism $\phi$ specified by the 9-tuple
$$(1, 11, 6, 3, 15, 9, 7, 9, 8),$$
with the endomorphism $\theta$ specified by the 9-tuple
$$(7, 11, 9, 1, 15, 8, 3, 9, 6).$$
Notice that these multisets are the same.  Thus, the sets of divergence 
rates are also identical.  The \NBC groups 
$\Gamma_\phi$ and $\Gamma_\theta$ cannot be distinguished via a comparison 
of divergence rates $\D$.

Next, we consider the growth spaces associated to these groups.  Denote by 
$\g_n$ the growth space of $\Gamma_\phi$ defined by $\lambda = 2^n$, and 
similarly denote by 
$\g'_n$ the growth space of $\Gamma_\theta$ defined by $\lambda = 2^n$.  We 
denote the Lie algebra of the Heisenberg group by $\h$, and the Lie algebra 
of $\R^n$ by $\re^n$.  Then,

\begin{align*}
\g_{1} &= \langle a_1  \rangle \cong \re&
\g'_{1} &= \langle a_4 \rangle \cong \re\\
\g_{3} &= \langle  a_1, a_4\rangle \cong \re^2&
\g'_{3} &= \langle  a_4, a_7\rangle \cong \re^2\\
\g_{6} &= \langle  a_1, a_4, a_3 \rangle \cong \re^3&
\g'_{6} &= \langle  a_4, a_7, a_9 \rangle \cong \re^3\\
\g_{7} &= \langle  a_1, a_4, a_3, a_7\rangle \cong \re^4&
\g'_{7} &= \langle  a_4, a_7, a_9, a_1 \rangle \cong \re^4\\
\g_{8} &= \langle  a_1, a_4, a_3, a_7, a_9\rangle \cong \re^5&
\g'_{8} &= \langle    a_4, a_7, a_9, a_1, a_6 \rangle \cong \re^5\\
\g_{9} &= \langle a_1, a_4, a_3, a_7, a_9, a_6, a_8  \rangle \cong \h 
\times \re^3&
\g'_{9} &= \langle   a_4, a_7, a_9, a_1, a_6, a_3, a_8\rangle \cong \h 
\times \re^3\\
\g_{11} &= \langle a_1, a_4, a_3, a_7, a_9, a_6, a_8, a_2 \rangle \cong 
\h^2 \times \re&
\g'_{11} &= \langle a_4, a_7, a_9, a_1, a_6, a_3, a_8, a_2 \rangle \cong 
\h^2 \times \re\\
\g_{15} &\cong \h^3&
\g'_{15} &\cong \h^3\\
\end{align*}

Notice that, as we consider the increasing sequence 
of growth spaces, a copy of $\h$ is introduced to a growth space whenever we 
reach the largest of the three eigenvalues associated to the terms of a 
component group $H$.  For both groups, this occurs at
at $\lambda = 2^9, 2^{11}$, and $2^{15}$.

These growth spaces are all isomorphic.  Therefore, they do not tell us that 
these groups are not \qic.  In fact, these groups have the same permuted 
absolute Jordan form, and thus are \qic.  Yet, each of the groups is a product 
of 3 \NBC groups, one of the six of which are pairwise \qic.

\subsection{Another Example}
\label{subsection:example2}
Recall the four-step nilpotent group $G$ described at the end of Subsection
\ref{subsection:kar} by Equation \ref{eqn:4step}.  Compare the \NBC groups 
defined by the two endomorphisms $\phi$ and $\theta$ of $G$.  As in the 
previous subsection, each
endomophism acts by raising each generator to a power, which is itself a
power of two.  Thus, the endomorphism is determined by the growth rates of
each generator.  This can be calculated for a generator $g$ to be
$n = \frac{1}{i} \log_2 k$, where $\phi(g) = g^k$ and $g \in V_i$.
For each of the endomorphisms considered here, we have $n_g = 3$ for $g =
p,q,r,s,t$.  Thus each
endomorphism is specified by the 6-tuple $(n_x, n_y, n_z, n_a, n_b,
n_c)$.  Compare the endomorphism $\phi$ specified by the 6-tuple
$$(1, 5, 3, 2, 4, 3),$$
with the endomorphism $\theta$ specified by the 6-tuple
$$(1, 3, 2, 3, 5, 4).$$

Because these multisets are the same, the divergence rates $\D_\phi$ and
$\D_\theta$ also agree.  The long but straightfoward computation of growth
spaces shows that these also agree.  However, the matrices associated to
these groups do not have the same permuted absolute forms.  In fact, 
these groups are noticably different.  Although the growth rates are the 
same, they are associated with points at different levels of the group.

A proof that these groups are not \QIc would be further evidence for the 
conjecture that permuted absolute Jordan form is a \QI invariant.

\section{The Rigidity of Nilpotent-by-Cyclic
Groups}\label{section:rigidity} %Section 4

In this section, we prove Theorem \ref{theorem:rigidity}. 
The proof of the corresponding rigidity theorem for \ABC groups (Theorem 
1.2 in 
\cite{ABC}) proceeds in six
steps.  Steps 1-3 apply here directly. Step 4 is modified slightly, and Steps
5-6 are not applicable to the nilpotent case.  As a result, we have the
weaker conclusion that $G$ is (virtually-nilpotent)-by-cyclic, and not the
stronger commensurability result found in the \ABC case.

\begin{proof}[Proof of Theorem \ref{theorem:rigidity}]
\hfill

\noindent\textbf{Step 1.}
The action of $G$ on itself by left multiplication can be conjugated by
the quasi-isometry $G \to X_{N, \phi}$ to give a proper, cobounded
quasi-action of $G$ on $X_{N, \phi}$ (see \cite{FarbMosher:BSTwo},
Proposition 2.1). Since $[N: \phi(N)] > 1$ we may apply Theorem
\ref{theorem:horizontal} to conclude that the quasi-action of $G$ on
$X_{N,\phi}$ coarsely respects the fibers of the uniform metric fibration
$X_{N, \phi} \to T_{N, \phi}$.

\medskip
\noindent
\textbf{Step 2.}
Now we use the following result of Mosher, Sageev, 
and Whyte:

\begin{theorem}[\cite{MSW}, Theorem 1] 
Fix an integer $n \geq 0$ and let $\Gamma$ be a finite graph of coarse 
$PD(n)$ groups with bushy Bass-Serre tree.  Let $H$ be a 
finitely-generated group \QIc to $\pi_1\Gamma$.  Then $H$ is the fundamental group 
of a graph of groups with bushy Bass-Serre tree and with vertex and edge 
groups \QIc to those of $\Gamma$.
\end{theorem}

By Step 1, this result applies to the quasi-action of $G$ on $X_{N,
\phi}$, because $G$ is quasi-isometric to the finitely-presented group
$\Gamma_{N, \phi}$ and so $G$ is finitely presented. The fibers of the map
$X_{N, \phi} \to T_{N,\phi}$ are isometric to $N$, and it follows that $G$
is the fundamental group of a graph of groups with each vertex and edge
group quasi-isometric to a nilpotent group.

\medskip
\noindent
\textbf{Step 3.}
By Gromov's polynomial growth theorem
\cite{G}, any finitely-generated group
quasi-isometric to a nilpotent group is \emph{virtually nilpotent}; that
is, it has a finite index nilpotent subgroup. 
Thus $G$ is the fundamental group of a graph of groups
whose vertex and edge groups are virtually nilpotent.

\medskip
\noindent
\textbf{Step 4.}
Amenability is a quasi-isometry invariant, so since all
\NBC groups are amenable, $G$ must also be amenable. No amenable group has
a nonabelian free subgroup, so $G$ has no free nonabelian subgroup.

Therefore the tree on which $G$ acts 
has $1$ `in' and $k$ `out' branches at each vertex.  (If there were two 
of each, then there would be two translation axes in the group 
action, and then by the Ping-Pong Lemma, it would have a nonabelian free 
subgroup).  
The fact that $G$ acts on a tree of this form implies that $G$ is the 
ascending HNN extension of some virtually nilpotent group $N'$; that is, 
$\Gamma$ is (virtually-nilpotent)-by-cyclic.
\end{proof}

\bigskip
This theorem is false if the condition that $\Gamma$ be finitely 
presented is weakened to finitely generated.  In fact, Dioubina 
\cite{D} has found examples of finitely-generated groups which are 
\QIc to \NBC groups, but not themselves (virtually-nilpotent)-by-cyclic. 

In particular, she shows that $\Z \wr \Z$ is quasi-isometric to $(Z \oplus
F) \wr \Z$, where $F$ is a finite nonsolvable group.  However, the former
group is solvable, while the latter is not.

\smallskip	
\smallskip	
\smallskip	
\emph{Acknowledgements}

I owe many thanks to Benson Farb for his outstanding mathematical
mentoring and faithful encouragement.  As advisor, he went far above and beyond the call of duty.  I also appreciate the mathematical
assistance of Chris Connell, John Franks, Lee Mosher, Kevin Whyte, and
Dave Witte.

Thanks to Michael Mihalik and the mathematics department at Vanderbilt
University for hosting me during the completion of this research.

\pagebreak
 
\end{document}